\documentclass[12pt]{article}

\usepackage{a4wide,epsfig}
\usepackage{amsfonts,amsmath}
\usepackage{graphicx}

\newtheorem{theorem}{Theorem}[section]

\newtheorem{corollary}[theorem]{Corollary}
\newtheorem{lemma}[theorem]{Lemma}

\newtheorem{remark}[theorem]{Remark}

\newcommand{\proof} [1]
   { \noindent {\bf Proof.} #1 \hfill\rule{0.5em}{1.2ex} \par\medskip}

\newcommand{\R}{\mathbb{R}}
\newcommand{\N}{\mathbb{N}}

\numberwithin{equation}{section} 

\begin{document}

\setcounter{page}{1}

\title{On a modified Hilbert transformation, \\
the discrete inf-sup condition, and error estimates}
\author{Richard~L\"oscher, Olaf~Steinbach, Marco~Zank}
\date{Institut f\"ur Angewandte Mathematik, TU Graz, \\[1mm]
Steyrergasse 30, 8010 Graz, Austria}

\maketitle

\begin{abstract}
  In this paper, we analyze the discrete inf-sup condition and related
  error estimates for a modified Hilbert transformation as used in the
  space-time discretization of time-dependent partial differential
  equations. It turns out that the stability constant depends linearly
  on the finite element mesh parameter, but in most cases, we can show
  optimal convergence. We present a series of numerical experiments
  which illustrate the theoretical findings.
\end{abstract}

\section{Introduction}
The Hilbert transformation ${\mathcal{H}}$, e.g.,
\cite{ButzerTrebels:1968, King:2009}, is a useful tool in the mathematical
analysis of parabolic evolution equations, e.g.,
\cite{AuscherEgertNystrom:2020, Fontes:1996}. The particular feature
is that the space-time bilinear form
$\langle \partial_t u , {\mathcal{H}} v \rangle$ is symmetric and
elliptic, and therefore allows for the use of standard arguments in the
numerical analysis of related space-time finite element methods~\cite{Devaud:2020, LarssonSchwab:2015} on the unbounded time interval $(0,\infty)$. While the Hilbert transformation
${\mathcal{H}}$ is defined as Cauchy principal value over ${\mathbb{R}}$, we have recently introduced a modified
Hilbert transformation ${\mathcal{H}}_T$ acting on a finite
time interval $(0,T)$ in \cite{SteinbachZank:2020, Zank:2020}. The relation between the modified transformation
${\mathcal{H}}_T$ and the classical one was recently given
in \cite{Ferrari:2023, SteinbachMissoni:2023}. In \cite{SteinbachZank:2020, Zank:2020},
we have analyzed a related space-time finite element method for the
numerical solution of a heat equation with zero initial and Dirichlet
boundary conditions as a model problem. Extensions include temporal $hp$ approaches and graded meshes in space \cite{PerugiaSchwabZank:2023}, space-time finite element methods
for the Maxwell's equations \cite{HauserZank:2023}, and the design of
efficient direct solution methods \cite{LangerZank:2021}. While in
the continuous case
we were able to establish an inf-sup stability condition as an ingredient
to ensure unique solvability of the space-time variational formulation,
the derivation of error estimates for the space-time finite element
approximation on tensor-product meshes used the ellipticity estimate
for the temporal part only. Indeed, numerical results in
\cite[Remark~3.4.29]{Zank:2020} indicate that there does not hold a discrete inf-sup
stability condition which is uniform in the space-time finite
element mesh size $h$, although we have observed optimal orders
of convergence.

When considering the formulation of boundary integral equations for the
wave equation in one space dimension, the composition of the modified
Hilbert transformation ${\mathcal{H}}_T$ with the acoustic single layer
boundary integral operator $V$ becomes self-adjoint and elliptic, see 
\cite{SteinbachUrzuaZank:2022}. This was
the motivation to consider a space-time finite element method for the
wave equation using the modified Hilbert transformation
${\mathcal{H}}_T$ applied to the test function. While we could ensure
unique solvability for any choice of conforming space-time
finite element spaces, at that time we were not able to prove
convergence, although we have observed optimal rates of convergence
in our numerical
examples \cite{LoescherSteinbachZank:2022}. Again, and as in the parabolic
case, numerical results indicate, that the related discrete inf-sup
stability condition does not hold uniformly in the space-time finite
element mesh size.

The question of analyzing the finite element error in the case when the
discrete inf-sup stability condition does not hold uniformly in the finite
element mesh size is well studied in the literature, see, e.g., the recent
work \cite{BertrandBoffi:2022}, and the references given therein. In
engineering, this is known as patch test, e.g.,
\cite{BabuskaNarasimhan:1997,ZienkiewiczTaylor:1997}, but also some
limitations are well documented, e.g., \cite{Stummel:1979, Stummel:1980}.

In this paper, we consider a modified Hilbert transformation based projection,
i.e., a Galerkin projection on piecewise polynomials where the test
functions involve the modified Hilbert transformation~$\mathcal H_T$.
This requires the analysis of the temporal bilinear
form~$\langle u , {\mathcal{H}}_T v \rangle_{L^2(0,T)}$ which is non-negative for
$u=v$, and which is a necessary
ingredient in the numerical analysis of space-time
variational formulations for parabolic and hyperbolic evolution equations.
Numerical results indicate that the constant of the related inf-sup
condition is proportional to the finite element mesh size $h$, independent
of the polynomial degree of the finite element space. On the other hand,
in almost all numerical examples, we observe convergence rates as expected
from the approximation properties of the finite element spaces. But in the
case of a function with a singularity at the origin $t=0$, we observe a
convergence rate that is much less than expected. For ease of presentation,
in this paper, we provide a detailed analysis in the case of piecewise
constant basis functions only. But this approach can be extended to higher-order polynomial basis functions. In fact, we prove that the discrete
inf-sup constant is proportional to the finite element mesh size, as
already observed in the numerical examples. Moreover, we are able
to characterize the discrete inf-sup constant in terms of the Fourier
coefficients of finite element functions with respect to the generating
functions of the modified Hilbert transformation ${\mathcal{H}}_T$. In
a second step, we analyze the projection error $u-u_h$ with respect to
the error $u-Q_hu$ of the standard $L^2$ projection $Q_hu$. When assuming
some regularity on $u$, we prove some super-convergence for $u_h-Q_hu$.
This is due to an appropriate splitting of $u_h-Q_hu$,
and mapping properties of $Q_h {\mathcal{H}}_T^{-1}$. In particular for
$ u \in H^2(0,T)$ satisfying $\partial_t u(0)=0$, we prove optimal
convergence although the discrete inf-sup condition is mesh dependent.
In all other cases, we provide a detailed analysis to theoretically explain all
convergence results as observed in the numerical examples.

The rest of this paper is organized as follows:
In Section \ref{Section:Hilbert}, we recall the definition of the modified
Hilbert transformation ${\mathcal{H}}_T$ and its properties. In particular,
in Lemma~\ref{Lemma H positiv}
we prove that $\langle v , {\mathcal{H}}_T v \rangle_{L^2(0,T)} > 0$
whenever $0 \neq v \in H^s_{0,}(0,T)$ for some $s \in (0,1]$, i.e., $v$ is more
regular than just in $L^2(0,T)$. In addition, we also provide an alternative
proof for the relation between the modified Hilbert transformation
${\mathcal{H}}_T$ and the Hilbert transformation~${\mathcal{H}}$ which is
different from what was presented in \cite{Ferrari:2023}. The modified
Hilbert transformation based projection is introduced in Section
\ref{Section:Projection}, see \eqref{Hilbert L2 Projektion}, and we
discuss the application of more standard stability and error estimates.
Numerical examples indicate that the discrete inf-sup condition depends
linearly on the finite element mesh size $h$, but the error behaves
with optimal order in most cases. While these numerical experiments
are done for basis functions up to second order, all further considerations
are done for piecewise constant basis functions only. In Section~\ref{Section:inf-sup}, we provide a proof for the discrete inf-sup condition
with a mesh dependent stability constant, see Theorem~\ref{Theorem inf-sup}. The stability constant as given in \eqref{csuh}
includes a dependency on the Fourier coefficients of the finite element
function $u_h$ which explains the different behavior as observed in the
numerical experiments. Related error estimates are then derived in
Section~\ref{Section:Error}. With these results, we are able to explain all the
convergence results as observed in the numerical experiments. In Section~\ref{Section:Con}, we finish
with some conclusions and comments on ongoing work.
For completeness, we provide all the details of the more technical
computations in the appendix.

\section{A modified Hilbert transformation}\label{Section:Hilbert}
Let $T>0$ be a given time horizon. For a function $v \in L^2(0,T)$,
we consider the Fourier series
\begin{equation*}
v(t) = \sum\limits_{k=0}^\infty v_k
\sin \left( \left( \frac{\pi}{2} + k\pi \right) \frac{t}{T} \right), \quad
t \in (0,T),
\end{equation*}
where the Fourier coefficients are given by
\begin{equation} \label{Fourier Koeff}
v_k = \frac{2}{T} \int_0^T v(t) 
\sin \left( \left( \frac{\pi}{2} + k\pi \right) \frac{t}{T} \right) dt .
\end{equation}
With this, we introduce the modified Hilbert transformation
${\mathcal{H}}_T\colon \, L^2(0,T) \to L^2(0,T)$ as, see \cite{SteinbachZank:2020},
\begin{equation*}
({\mathcal{H}}_Tv)(t) := \sum\limits_{k=0}^\infty v_k
\cos \left( \left( \frac{\pi}{2} + k\pi \right) \frac{t}{T} \right), \quad t \in (0,T).
\end{equation*}
Note that by Parseval's theorem, we have
\begin{equation}\label{Parseval}
  \| {\mathcal H}_T v \|_{L^2(0,T)}^2 =
  \| v \|_{L^2(0,T)}^2 = \frac{T}{2} \, \sum\limits_{k=0}^\infty v_k^2 \, .
\end{equation}
The inverse of the modified Hilbert transformation ${\mathcal{H}}_T^{-1} \colon \, L^2(0,T) \to L^2(0,T)$ is given by
\begin{equation*}
({\mathcal{H}}_T^{-1}w)(t) = \sum\limits_{k=0}^\infty \overline{w}_k
\sin \left( \left( \frac{\pi}{2} + k\pi \right) \frac{t}{T} \right), \quad t \in (0,T),
\end{equation*}
where
\begin{equation}\label{Fourier coefficients cos}
\overline{w}_k = \frac{2}{T} \int_0^T w(t) 
\cos \left( \left( \frac{\pi}{2} + k\pi \right) \frac{t}{T} \right) dt 
\end{equation}
are the coefficients of the Fourier series
\begin{equation}\label{Fourier series cos}
w(t) = \sum\limits_{k=0}^\infty \overline{w}_k
\cos \left( \left( \frac{\pi}{2} + k\pi \right) \frac{t}{T} \right), \quad t \in (0,T),
\end{equation}
of a function $w \in L^2(0,T)$. As in \cite[Lemma 2.4]{SteinbachZank:2020}, we have
\begin{equation}\label{HT-HTinverse}
  \langle {\mathcal{H}}_T u , w \rangle_{L^2(0,T)} =
  \langle u , {\mathcal{H}}_T^{-1} w \rangle_{L^2(0,T)} \quad
  \mbox{for all} \; u, w \in L^2(0,T).
\end{equation}
In fact, ${\mathcal{H}}_T \colon \, L^2(0,T) \to L^2(0,T)$ is an isometry, which implies 
the inf-sup stability condition
\begin{equation*}
\| u \|_{L^2(0,T)} = \sup\limits_{0 \neq v \in L^2(0,T)}
\frac{\langle u , {\mathcal{H}}_T v \rangle_{L^2(0,T)}}{\| v \|_{L^2(0,T)}}
\quad \mbox{for all} \; u \in L^2(0,T) .
\end{equation*}
We denote by $H^s(0,T)$ for $s>0$ the usual Sobolev spaces with norm $\| \cdot \|_{H^s(0,T)}$. Further, we define the closed subspaces
\[
  H^1_{0,}(0,T) := \{ v \in H^1(0,T) : v(0)=0 \}, \quad
  H^1_{,0}(0,T) := \{ v \in H^1(0,T) : v(T)=0 \}
\]
of $H^1(0,T)$ endowed with the Hilbertian norms $\|\cdot\|_{H^1_{0,}(0,T)} := \|\cdot\|_{H^1_{,0}(0,T)}  := \|\partial_t (\cdot) \|_{L^2(0,T)}$ and, we introduce the interpolation spaces
\[
  H^s_{0,}(0,T) := [H^1_{0,}(0,T),L^2(0,T)]_s, \quad
  H^s_{,0}(0,T) := [H^1_{,0}(0,T),L^2(0,T)]_s
\]
for $s \in (0,1)$, see \cite[Section~2.2]{Zank:2020} for further references.  
We equip these interpolation spaces with the norms in Fourier representation, i.e.,
\begin{align*}
    \| v \|_{H^s_{0,}(0,T)} :=& \left( \frac{T}{2} \sum\limits_{k=0}^\infty 
      \left[ \frac{1}{T} \left( \frac{\pi}{2} + k\pi \right) \right]^{2s} v_k^2 \right)^{1/2}, \\
    \| w \|_{H^s_{,0}(0,T)} :=& \left( \frac{T}{2} \sum\limits_{k=0}^\infty 
      \left[ \frac{1}{T} \left( \frac{\pi}{2} + k\pi \right) \right]^{2s} \overline{w}_k^2 \right)^{1/2}
\end{align*}
for $v \in H^s_{0,}(0,T)$ with Fourier coefficients $v_k$ given in
\eqref{Fourier Koeff}, and $w \in H^s_{,0}(0,T)$ with Fourier
coefficients $\overline{w}_k$ given in \eqref{Fourier coefficients cos}.
Further, for $s \in (0,1)$, $\langle \cdot,\cdot \rangle_{(0,T)}$ denotes the duality pairing
in $H^s_{,0}(0,T)$ and $[H^s_{,0}(0,T)]'$ as continuous extension of
the $L^2(0,T)$ inner product $\langle \cdot,\cdot \rangle_{L^2(0,T)}$. Here, for $s \in (0,1)$,
the dual space $[H^s_{,0}(0,T)]'$ is endowed with the norm
\begin{equation*}
  \| z \|_{[H^s_{,0}(0,T)]'} = \sup_{0 \neq w \in H^s_{,0}(0,T)}
  \frac{\langle z, w \rangle_{(0,T)}}{ \| w \|_{H^s_{,0}(0,T)} },
  \quad z \in [H^s_{,0}(0,T)]'.
\end{equation*}
With this notation, we have the ellipticity
\[
  \langle \partial_t v , {\mathcal{H}}_T v \rangle_{(0,T)}  \, = \,
  \| v \|^2_{H^{1/2}_{0,}(0,T)} 
\quad \mbox{for all} \; v \in H^{1/2}_{0,}(0,T),
\]
see \cite[Equation~(2.9)]{SteinbachZank:2020}. Further, for $s \in [0,1]$, the modified Hilbert transformation is an isometry as mapping $\mathcal H_T \colon \, H^s_{0,}(0,T) \to H^s_{,0}(0,T)$, satisfying
\begin{equation*}
  \| \mathcal H_T v \|_{H^s_{,0}(0,T)} = \| v \|_{H^s_{0,}(0,T)} \quad \text{for all } v \in H^s_{0,}(0,T),
\end{equation*}
which implies a second inf-sup stability condition,
\begin{equation*}
\| u \|_{[H^{1/2}_{,0}(0,T)]'} =
\sup\limits_{0 \neq v \in H^{1/2}_{0,}(0,T)}
\frac{\langle u , {\mathcal{H}}_T v \rangle_{(0,T)}}
{\| v \|_{H^{1/2}_{0,}(0,T)}} \quad
\mbox{for all} \; u \in [H^{1/2}_{,0}(0,T)]' .
\end{equation*}

\begin{remark} \label{Bem:infsupH12}
  Note that for $0 \neq u \in L^2(0,T)$, the function
  $\overline{v}(t) := \int_0^t u(s) \, ds$, $t \in (0,T)$, yields
\begin{align*}
\langle u , {\mathcal{H}}_T \overline{v} \rangle_{(0,T)} = 
\langle \partial_t \overline{v} , {\mathcal{H}}_T \overline{v} \rangle_{(0,T)}
 = 
\| \overline{v} \|^2_{H^{1/2}_{0,}(0,T)}  &= 
\| \overline{v} \|_{H^{1/2}_{0,}(0,T)} 
\| \partial_t \overline{v} \|_{[H^{1/2}_{,0}(0,T)]'} \\
&= 
\| \overline{v} \|_{H^{1/2}_{0,}(0,T)} 
\| u \|_{[H^{1/2}_{,0}(0,T)]'}
\end{align*}
and thus, 
\begin{equation*}
\| u \|_{[H^{1/2}_{,0}(0,T)]'} = 
\frac{\langle u , {\mathcal{H}}_T \overline{v} \rangle_{(0,T)}}
{\| \overline{v} \|_{H^{1/2}_{0,}(0,T)}} =
\sup\limits_{0 \neq v \in H^{1/2}_{0,}(0,T)}
\frac{\langle u , {\mathcal{H}}_T v \rangle_{(0,T)}}
{\| v \|_{H^{1/2}_{0,}(0,T)}}.
\end{equation*}
\end{remark}

\noindent
In \cite[Lemma 2.6]{SteinbachZank:2020}, it was shown that
\[
  \langle v , {\mathcal{H}}_T v \rangle_{L^2(0,T)} \geq 0 \quad
  \mbox{for all} \; v \in L^2(0,T).
\]
In fact, for a bit more regular functions, we have the following result:

\begin{lemma}\label{Lemma H positiv}
  For $0 \neq v \in H^s_{0,}(0,T) =  [H^1_{0,}(0,T),L^2(0,T)]_s$ with
  $s \in (0,1]$, the inequality
  \[
\langle v , {\mathcal{H}}_T v \rangle_{L^2(0,T)} > 0
  \]
  holds true.
\end{lemma}

\proof{By using the representations
\[
v(t) = \sum\limits_{k=0}^\infty v_k 
\sin \left( \Big( \frac{\pi}{2} + k\pi \Big) \frac{t}{T} \right), \quad
({\mathcal{H}}_Tv)(t) = \sum\limits_{\ell=0}^\infty 
v_\ell \cos \left( \Big( \frac{\pi}{2} + \ell\pi \Big) \frac{t}{T} \right),
\]
it follows from the proof of \cite[Lemma~2.6]{SteinbachZank:2020} that
\begin{equation}
	\langle v , {\mathcal{H}}_T v \rangle_{L^2(0,T)} = \frac{T}{\pi} \lim_{N \to \infty}  \left[ \int_0^1 \left( \sum\limits_{i=0}^N v_{2i} x^{2i} \right)^2 dx + \int_0^1 \left( \sum\limits_{i=0}^N v_{2i+1} x^{2i+1} \right)^2 dx \right]. \label{PosdefGrenzwertIntegral}
\end{equation}
To interchange the limit processes and integral signs in \eqref{PosdefGrenzwertIntegral}, the \textit{Theorem of Lebesgue} is applied. For the first part in \eqref{PosdefGrenzwertIntegral}, define the functions
\begin{equation*}
    f_{1,N} \colon \, (0,1) \to [0,\infty], \quad f_{1,N}(x) := \left( \sum\limits_{i=0}^N v_{2i} x^{2i} \right)^2 \quad \text{for } N \in \N_0,
\end{equation*}
and
\begin{equation*}
    f_1\colon \, (0,1) \to [0,\infty], \quad f_1(x) := \left( \sum\limits_{i=0}^\infty v_{2i} x^{2i} \right)^2.
\end{equation*}
The Cauchy--Schwarz inequality yields
\begin{equation*}
    f_{1,N}(x) = | f_{1,N}(x) | \leq \frac{T^{1-2s}}{2} \sum\limits_{i=0}^\infty v_{2i}^2 \left(\frac{\pi}{2} + 2 i \pi \right)^{2 s}   \frac{2}{T^{1-2s}}  \sum\limits_{i=0}^\infty \frac{x^{4i}}{ \left(\frac{\pi}{2} + 2i \pi \right)^{2 s}  } \leq \| v \|_{H^s_{0,}(0,T)}^2 \cdot g_1(x)
\end{equation*}
for all $x \in (0,1)$ and for all $N \in \N_0$, where
\begin{equation*}
    g_1\colon \, [-1,1] \to [0,\infty], \quad g_1(x) := \frac{2}{T^{1-2s}} \sum\limits_{i=0}^\infty \frac{x^{4i}}{ \left(\frac{\pi}{2} + 2i \pi \right)^{2 s}  } = \frac{2^{2s+1}}{T^{1-2s} \pi^{2s}} \sum\limits_{i=0}^\infty \frac{x^{4i}}{ \left(1 + 4i \right)^{2 s}  }.
\end{equation*}
For $x\in (-1,1)$, the estimates
\begin{equation*}
 0 \leq \sum\limits_{i=0}^\infty \frac{x^{4i}}{ \left(1 + 4i \right)^{2 s}  } \leq \sum\limits_{i=0}^\infty (x^4)^i = \frac{1}{1-x^4}
\end{equation*}
hold true, i.e., the power series $\sum\limits_{i=0}^\infty \frac{x^{4i}}{ \left(1 + 4i \right)^{2 s}  }$ converges absolutely for $x \in (-1,1)$. Hence, for any $y \in (0,1)$, termwise integration is applicable and gives
\begin{equation*}
  \int_0^y | g_1(x)| dx  = \frac{2^{2s+1}}{T^{1-2s} \pi^{2s}} \sum\limits_{i=0}^\infty \int_0^y \frac{x^{4i}}{ \left(1 + 4i \right)^{2 s}  } dx = \frac{2^{2s+1}}{T^{1-2s} \pi^{2s}} \sum\limits_{i=0}^\infty \frac{y^{4i+1}}{ \left(1 + 4i \right)^{2 s+1}  },
\end{equation*}
where all occurring series converge absolutely. For the endpoint $y=1$, we have that the series $\sum\limits_{i=0}^\infty \frac{1}{ \left(1 + 4i \right)^{2 s+1}  }$ converges if and only if $s>0$, since it can be estimated by general harmonic series. Hence, \textit{Abel's Theorem} is applicable, which states
\begin{equation*}
  \int_0^1 | g_1(x) | dx = \lim_{y \to 1^-}\int_0^y | g_1(x) | dx  = \frac{2^{2s+1}}{T^{1-2s} \pi^{2s}} \sum\limits_{i=0}^\infty \frac{1}{ \left(1 + 4i \right)^{2 s+1}  }.
\end{equation*}
In other words, the function $g_1$ is integrable on the interval $(0,1)$ if and only if $s > 0.$ Analogous results hold true for the second part in \eqref{PosdefGrenzwertIntegral} with the involved functions
\begin{align*}
    &f_{2,N} \colon \, (0,1) \to [0,\infty], \quad f_{2,N}(x) := \left( \sum\limits_{i=0}^N v_{2i+1} x^{2i+1} \right)^2 \quad \text{for } N \in \N_0, \\
    &f_2\colon \, (0,1) \to [0,\infty], \quad f_2(x) := \left( \sum\limits_{i=0}^\infty v_{2i+1} x^{2i+1} \right)^2
\end{align*}
and the function $g_2\colon \, [-1,1] \to [0,\infty],$
\begin{equation*}
    g_2(x) := \frac{2}{T^{1-2s}} \sum\limits_{i=0}^\infty \frac{x^{4i+2}}{ \left(\frac{\pi}{2} + (2i+1) \pi \right)^{2 s}  } = \frac{2^{2s+1}}{T^{1-2s} \pi^{2s}} \sum\limits_{i=0}^\infty \frac{x^{4i+2}}{ \left(1 + (4i +2) \right)^{2 s}  },
\end{equation*}
which is integrable on $(0,1)$ if and only if $s>0.$ Summarizing, all assumptions of the \textit{Theorem of Lebesgue} are satisfied for \eqref{PosdefGrenzwertIntegral}. Thus, interchanging the limit processes and integral signs in \eqref{PosdefGrenzwertIntegral} gives
\begin{equation*}
	\langle v , {\mathcal{H}}_T v \rangle_{L^2(0,T)} = \frac{T}{\pi}  \left[ \int_0^1 \left( \sum\limits_{i=0}^\infty v_{2i} x^{2i} \right)^2 dx + \int_0^1 \left( \sum\limits_{i=0}^\infty v_{2i+1} x^{2i+1} \right)^2 dx \right].
\end{equation*}
The last line is strictly positive, since otherwise, the relation
\begin{equation*}
 \forall x \in (0,1) \colon \quad  \sum\limits_{i=0}^\infty v_{2i} x^{2i} = \sum\limits_{i=0}^\infty v_{2i+1} x^{2i+1} = 0
\end{equation*}
would lead, by the \textit{Identity Theorem of Power Series}, to $v_k = 0$ for all $k \in \N_0$, i.e., a contradiction to $v \neq 0$.
}

\noindent
For $v \in L^2(0,T)$, the closed representation
\cite[Lemma 2.1]{SteinbachZank:2021}
\[
  ({\mathcal{H}}_Tv)(t) = \frac{1}{2T} \, \mathrm{v.p.}
  \int_0^T \left[
    \frac{1}{\sin \left( \frac{\pi}{2} \frac{s-t}{T}\right)}
    +
    \frac{1}{\sin \left( \frac{\pi}{2} \frac{s-t}{T}\right)}
  \right] \, v(s) \, ds \quad \mbox{for} \; t \in (0,T)
\]
as Cauchy principal value holds true, where additional integral representations are contained in \cite{Zank:2023}.
Recall that the classical Hilbert transformation \cite{ButzerTrebels:1968}
is given as
\[
  ({\mathcal{H}}\varphi)(t) = \frac{1}{\pi} \, \mathrm{v.p.}
  \int_\R \frac{\varphi(s)}{t-s} \, ds \quad
  \mbox{for} \; t \in \R .
\]
In \cite[Lemma 4.1]{SteinbachMissoni:2023}, it was shown that the
modified Hilbert transformation ${\mathcal{H}}_T$ differs from the classical
Hilbert transformation ${\mathcal{H}}$ by a compact perturbation
$B \colon \, L^2(0,T) \to H^1(0,T)$, i.e.,
\[
  ({\mathcal{H}}_T\varphi)(t) = -
  ({\mathcal{H}}\overline{\varphi})(t) +
  (B\varphi)(t), \quad t \in (0,T),
\]
where
\[
  \overline{\varphi}(s) :=
  \left \{
    \begin{array}{lcl}
      \varphi(s) & & \mbox{for} \; s \in (0,T), \\[1mm]
      \varphi(2T-s) & & \mbox{for} \; s \in (T,2T), \\[1mm]
      - \varphi(-s) & & \mbox{for} \; s \in (-T,0), \\[1mm]
      - \varphi(2T+s) & & \mbox{for} \; s \in (-2T,-T), \\[1mm]
      0 & & \mbox{otherwise}
    \end{array} \right.
\]
is a double reflection of the given function $\varphi \in L^2(0,T)$.
As it was recently shown in \cite{Ferrari:2023}, the modified Hilbert
transformation ${\mathcal{H}}_T$ coincides with the classical Hilbert
transformation~${\mathcal{H}}$, when the extension of $\varphi \in L^2(0,T)$
is defined accordingly. While there are already three different proofs
of this property in \cite{Ferrari:2023}, here, we present a different
way of proof. To this end, we introduce some useful notation.
Namely, we denote for $k \in {\mathbb{Z}}$ by 
\[
\mathcal{E}_{o} \varphi(t) := \left\{\begin{array}{lcl}
	\varphi(s-4kT)&& \mbox{for}\; s\in [4kT,(4k+1)T),\\[1mm]
	\varphi((4k+2)T-s) && \mbox{for}\; s\in[(4k+1)T,(4k+2)T),\\[1mm]
	-\varphi(4kT-s) && \mbox{for} \; s\in[(4k-1)T,4kT),\\[1mm]
	-\varphi(s-(4k-2)T) && \mbox{for}\; s\in[(4k-2)T,(4k-1)T)
\end{array} \right.
\]
the odd extension of a given function $\varphi\in L^2(0,T)$ to a function on $\R$. In particular,
note that 
\begin{equation}\label{eq:odd-extension-sin}
  \mathcal{E}_o \sin\left(\Big(\frac{\pi}{2}+k\pi\Big)\frac{\cdot}{T} \right)(t) =
  \sin\left(\Big(\frac{\pi}{2}+k\pi\Big)\frac{t}{T} \right), \quad
  t\in\mathbb{R}.
\end{equation} 
Moreover, denote the restriction operator by
$\mathcal{R}_{(0,T)}\widetilde{\varphi}(t) = \widetilde{\varphi}_{|(0,T)}(t)$,
$t\in (0,T)$, for a function $\widetilde{\varphi} \colon \, \R \to \R$. Then, the modified Hilbert transformation admits the
representation
\begin{equation*}
  \mathcal{H}_T\varphi =
  -\mathcal{R}_{(0,T)}\mathcal{H}\mathcal{E}_o\varphi \quad \text{for } \varphi \in L^2(0,T).
\end{equation*}
To prove this, recall that each $\varphi \in L^2(0,T)$ admits the
representation 
\begin{equation*}
  \varphi(t) = \sum_{k=0}^\infty \varphi_k
  \sin\left(\Big(\frac{\pi}{2}+k\pi\Big)\frac{t}{T}\right),\quad
  \varphi_k = \frac{2}{T}\int _0^T \varphi(t)\sin\left(
    \Big(\frac{\pi}{2}+k\pi\Big)\frac{t}{T}\right)\, dt.
\end{equation*} 
It is well-known that $(\mathcal{H}\sin(a\cdot))(t) = -\cos(at)$ for $t \in \R$ and $a>0$,
see \cite[(3.110), p. 103]{King:2009}. Using this, and \eqref{eq:odd-extension-sin}, we compute that 
\begin{align*}
  (\mathcal{H}\mathcal{E}_o\varphi)(t)
  &=\Big(\mathcal{H}\mathcal{E}_o\sum_{k=0}^\infty \varphi_k
        \sin\left(\Big(\frac{\pi}{2}+k\pi\Big)\frac{\cdot}{T}\right)(t)
        \, = \,
	\sum_{k=0}^\infty \varphi_k \Big(\mathcal{H}\mathcal{E}_o
        \sin\left(\Big(\frac{\pi}{2}+k\pi\Big)\frac{\cdot}{T}\right)(t) \\
  &=\sum_{k=0}^\infty \varphi_k \Big(\mathcal{H}
        \sin\left(\Big(\frac{\pi}{2}+k\pi\Big)\frac{\cdot}{T}\right)(t)
        \, = \,
	-\sum_{k=0}^\infty \varphi_k \cos\left(
        \Big(\frac{\pi}{2}+k\pi\Big)\frac{t}{T}\right)\\
	&=-(\mathcal{H}_T\varphi)(t)
\end{align*}
for $t \in (0,T)$. Analogously, using the even extension, for $k \in {\mathbb{Z}}$, 
\[
  \mathcal{E}_{e} \varphi(t) := \left\{\begin{array}{lcl}
	\varphi(s-4kT)&& \mbox{for}\; s\in [4kT,(4k+1)T),\\[1mm]
	-\varphi((4k+2)T-s) && \mbox{for}\; s\in[(4k+1)T,(4k+2)T),\\[1mm]
	\varphi(4kT-s) && \mbox{for} \; s\in[(4k-1)T,4kT),\\[1mm]
	-\varphi(s-(4k-2)T) && \mbox{for}\; s\in[(4k-2)T,(4k-1)T)
\end{array} \right.
\]
for which 
\begin{equation*}
  \mathcal{E}_e \cos\left(\Big(\frac{\pi}{2}+k\pi\Big)\frac{\cdot}{T} \right)(t) =
  \cos\left(\Big(\frac{\pi}{2}+k\pi\Big)\frac{t}{T} \right), \quad
  t\in\mathbb{R},
\end{equation*}
together with the cosine expansion 
\begin{equation*}
  \varphi(t) = \sum_{k=0}^\infty \overline{\varphi}_k
  \cos\left(\Big(\frac{\pi}{2}+k\pi\Big)\frac{t}{T}\right),\quad
  \overline{\varphi}_k = \frac{2}{T}\int _0^T \varphi(t)
  \cos\left(\Big(\frac{\pi}{2}+k\pi\Big)\frac{t}{T}\right)\, dt,
\end{equation*}
and the property that $(\mathcal{H}\cos(a\cdot))(t) = \sin(at)$, $t \in \R$, $a>0$, we derive that 
\begin{equation*}
	\mathcal{H}_T^{-1} \varphi = \mathcal{R}_{(0,T)}\mathcal{H}\mathcal{E}_e\varphi \quad \text{for } \varphi \in L^2(0,T).
\end{equation*}
To check the consistency, note that the Hilbert transformation of even, periodic functions is odd with the same period and vice versa, see \cite[Section 4.2]{King:2009}. Thus, the equality
\begin{equation*}
	\mathcal{E}_e\mathcal{R}_{(0,T)}\mathcal{H}\mathcal{E}_o\varphi = \mathcal{H}\mathcal{E}_o\varphi \quad \text{for } \varphi \in L^2(0,T)
\end{equation*}
holds true. Using that $\mathcal{H}^2 = -\mathrm{id}$, it is easy to check that 
\begin{equation*}
	\mathcal{H}_T^{-1}\mathcal{H}_T\varphi = -\mathcal{R}_{(0,T)}\mathcal{H}\mathcal{E}_e\mathcal{R}_{(0,T)}\mathcal{H}\mathcal{E}_o\varphi = -\mathcal{R}_{(0,T)}\mathcal{H}\mathcal{H}\mathcal{E}_o\varphi = \varphi \quad \text{for } \varphi \in L^2(0,T).
\end{equation*}
The relation $\mathcal{H}_T\mathcal{H}_T^{-1}\varphi = \varphi$ for $\varphi \in L^2(0,T)$ follows in the same way.

\section{A modified Hilbert transformation based projection}
\label{Section:Projection}
For the finite time interval $(0,T)$ and a given discretization parameter
$n \in \N$, we consider a uniform decomposition of $(0,T)$ into $n$
finite elements $(t_{i-1},t_i)$ of mesh size $h=T/n$ with nodes
$t_i=ih$, $i=0,1,\ldots,n$. With respect to this mesh,
we introduce a conforming finite element space 
\begin{equation} \label{Vh}
V_h := S_h^\nu(0,T) = \mbox{span} \; \{ \psi^\nu_i \}_{i=1}^\mathrm{dof}\subset L^2(0,T)
\end{equation}
of either piecewise constant functions ($\nu=0$ and $\mathrm{dof} = n$) with basis
\[
  \psi^0_i(t) =
  \left \{
    \begin{array}{ccl}
      1 & & \mbox{for} \; t \in (t_{i-1},t_i), \\[1mm]
      0 & & \mbox{otherwise},
    \end{array} \right.     
\]
or piecewise linear continuous functions ($\nu=1$ and $\mathrm{dof} = n$) with basis
\[
  \psi_i^1(t) =
  \left \{
    \begin{array}{ccl}
      (t-t_{i-1})/h & & \mbox{for} \; t \in (t_{i-1},t_i], \\[1mm]
      (t_{i+1}-t)/h & & \mbox{for} \; t \in (t_i,t_{i+1}), \\[1mm]
      0 & & \mbox{otherwise}
    \end{array}
  \right.
\]
for $i=1,\dots,n-1$ and
\[
  \psi_n^1(t) =
  \left \{
    \begin{array}{ccl}
      (t-t_{n-1})/h & & \mbox{for} \; t \in (t_{n-1},t_n], \\[1mm]
      0 & & \mbox{otherwise},
    \end{array}
  \right.
\]
or piecewise quadratic continuous functions ($\nu=2$ and $\mathrm{dof} = 2n$) fulfilling homogeneous initial conditions. Thus, $S_h^\nu(0,T) \subset H^1_{0,}(0,T)$ for $\nu=1$ or $\nu=2$.

\begin{figure}[ht]
  \centering
  \begin{tabular}{cc}
      \includegraphics[width = 0.45\textwidth]{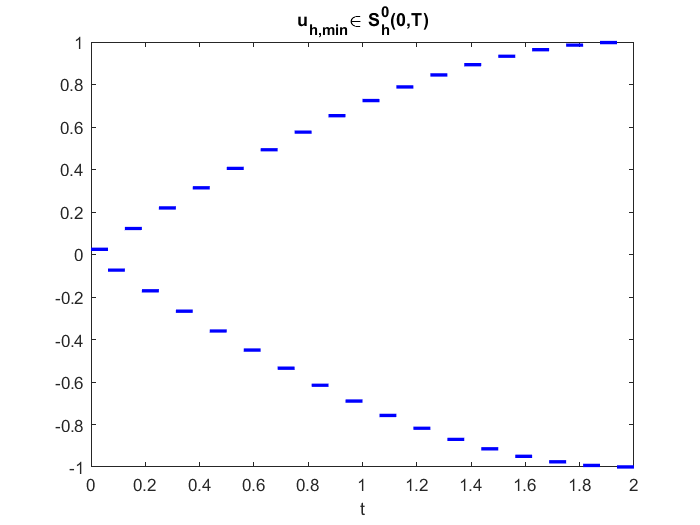}
      &
      \includegraphics[width = 0.45\textwidth]{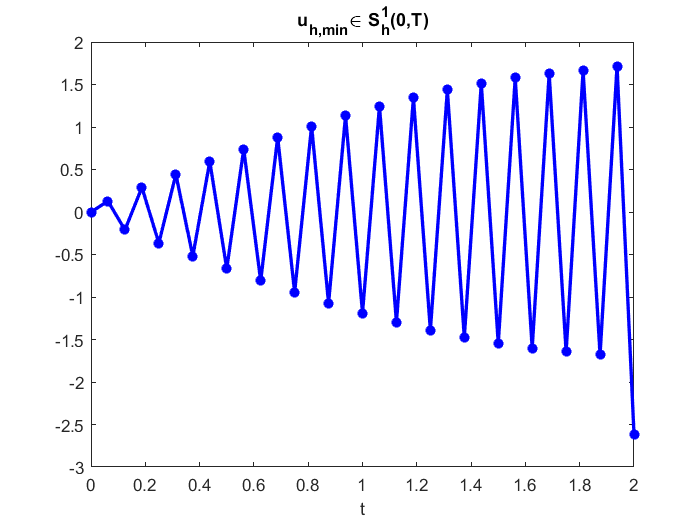}
  \end{tabular}
  \caption{Functions $u_{h,\text{min}}$ realizing \eqref{inf sup L2 diskret} with the smallest inf-sup constant $c_S$ for $N=32$ elements for $T=2$.} \label{Fig:infsup}
\end{figure}

Next, for given $u \in L^2(0,T)$, we consider the modified
Hilbert transformation based projection to find $u_h \in V_h$ such that
\begin{equation}\label{Hilbert L2 Projektion}
\langle u_h , {\mathcal{H}}_T v_h \rangle_{L^2(0,T)} =
\langle u , {\mathcal{H}}_T v_h \rangle_{L^2(0,T)} \quad 
\mbox{for all} \; v_h \in V_h .
\end{equation}
The variational formulation (\ref{Hilbert L2 Projektion}) is equivalent to 
a linear system of algebraic equations, 
\[
B_h^\nu \underline{u}=\underline{f}^\nu,
\]
where the Hilbert mass matrix $B_h^\nu$ is defined by its matrix entries
\[
B^\nu_h[j,i] = \langle \psi^\nu_i, {\mathcal{H}}_T \psi^\nu_j \rangle_{L^2(0,T)}
\quad \mbox{for} \; i,j=1,\ldots,\mathrm{dof},
\]
and the vector $\underline{f}^\nu$
of the right-hand side is given by its entries
\[
  f^\nu_j = \langle u , {\mathcal{H}}_T \psi^\nu_j \rangle_{L^2(0,T)} \quad
  \mbox{for} \; j=1,\ldots,\mathrm{dof}.
\]
Due to Lemma \ref{Lemma H positiv}, we have
\[
  \underline{v}^\top B_h^\nu \underline{v} =
  \langle v_h , {\mathcal{H}}_T v_h \rangle_{L^2(0,T)} > 0
  \quad \mbox{for all} \; \R^{\mathrm{dof}} \ni \underline{v}
  \leftrightarrow v_h \in V_h, v_h \neq 0,
\]
from which unique solvability of \eqref{Hilbert L2 Projektion} follows.
To prove related error estimates, we need to ensure 
the discrete inf-sup stability condition
\begin{equation}\label{inf sup L2 diskret}
c_S \, \| u_h \|_{L^2(0,T)} \leq \sup\limits_{0 \neq v_h \in V_h} 
\frac{\langle u_h , {\mathcal{H}}_Tv_h \rangle_{L^2(0,T)}}
{\| v_h \|_{L^2(0,T)}} \quad \mbox{for all} \; u_h \in V_h
\end{equation}
from which we also derive the a priori error estimate
\cite[Theorem~2]{XuZikatanov:2003}, i.e., Céa's lemma,
\begin{equation}\label{Cea}
\| u - u_h \|_{L^2(0,T)} \leq \frac{1}{c_S}\,
\inf\limits_{v_h \in V_h} \| u - v_h \|_{L^2(0,T)} \, .
\end{equation}
To motivate our theoretical considerations, let us first consider some 
numerical examples, where the calculation of the matrix $B^\nu_h$ and right-hand side $\underline{f}^\nu$ is done as proposed in \cite{Zank:2021}. In Table~\ref{Table inf sup constant L2}, we present
numerical results for the stability constant $c_S$ of the inf-sup stability
condition~\eqref{inf sup L2 diskret}, where for $T=2$,
the function $u_{h,\text{min}}$ realizing \eqref{inf sup L2 diskret} with the smallest inf-sup constant $c_S$ is depicted in Figure~\ref{Fig:infsup}. We observe that $c_S$ is mesh dependent.
In particular, we have $c_S \approx 0.426 \cdot h$ for $\nu=0$.

\begin{table}[ht]
\centering
\begin{tabular}{rlcccccc}
\hline
 & & \multicolumn{2}{c}{$\nu=0$} & \multicolumn{2}{c}{$\nu=1$} 
& \multicolumn{2}{c}{$\nu=2$} \\
\hline
$n$ & $h$ & $c_S$ & $c_S/h$ & $c_S$ & $c_S/h$ & $c_S$ & $c_S/h$ \\
\hline
   2 & 1.0      & 0.411711 & 0.412 & 0.515034 & 0.515 & 0.429033 & 0.429 \\
   4 & 0.5      & 0.211292 & 0.423 & 0.344142 & 0.688 & 0.271686 & 0.543 \\
   8 & 0.25     & 0.106338 & 0.425 & 0.204556 & 0.818 & 0.155494 & 0.622 \\
  16 & 0.125    & 0.053256 & 0.426 & 0.112324 & 0.899 & 0.083498 & 0.668 \\
  32 & 0.0625   & 0.026639 & 0.426 & 0.058935 & 0.943 & 0.043295 & 0.693 \\
  64 & 0.03125  & 0.013321 & 0.426 & 0.030192 & 0.966 & 0.022047 & 0.705 \\
 128 & 0.015625 & 0.006661 & 0.426 & 0.015281 & 0.978 & 0.011125 & 0.712 \\
 256 & 0.007812 & 0.003330 & 0.426 & 0.007687 & 0.984 & 0.005588 & 0.715 \\
 512 & 0.003906 & 0.001665 & 0.426 & 0.003855 & 0.987 & 0.002800 & 0.717 \\
1024 & 0.001953 & 0.000833 & 0.426 & 0.001931 & 0.988 & 0.001402 & 0.718 \\
2048 & 0.000977 & 0.000416 & 0.426 & 0.000966 & 0.989 & 0.000701 & 0.718 \\
\hline
\end{tabular}
\caption{Numerical results for the stability constant $c_S$ in
\eqref{inf sup L2 diskret} with $T=2$.}\label{Table inf sup constant L2}
\end{table}

From \eqref{Cea} and using the approximation properties of $V_h$ when
assuming $u \in H^s(0,T)$ for $s \in [0,\nu+1]$,
we then conclude the error estimate
\[
\| u - u_h \|_{L^2(0,T)} \leq c \, h^{s-1} \, \|u\|_{H^s(0,T)} \, ,
\]
in particular, we may lose one order in $h$. However, numerical results
indicate a rather different error behavior. As a first example, we
consider the regular function $u(t)=\sin \big( \frac{\pi}{4} t \big)$
for $t \in (0,T)$ with $T=2$, where we observe the optimal order of
convergence as expected from the approximation
properties of $V_h$, see Table~\ref{Table Fehler regular}.

\begin{table}[ht]
\centering
\begin{tabular}{rcccccc}
\hline
 & \multicolumn{2}{c}{$\nu=0$} 
 & \multicolumn{2}{c}{$\nu=1$} 
 & \multicolumn{2}{c}{$\nu=2$} \\ \hline
$n$ & $\| u - u_h \|_{L^2(0,T)}$ & eoc
    & $\| u - u_h \|_{L^2(0,T)}$ & eoc
    & $\| u - u_h \|_{L^2(0,T)}$ & eoc \\ \hline
   2 & 3.983 --1 &      & 2.605 --2 &      & 3.893 --3 & \\
   4 & 1.825 --1 & 1.13 & 5.933 --3 & 2.13 & 5.465 --4 & 2.83 \\
   8 & 8.971 --2 & 1.02 & 1.448 --3 & 2.03 & 7.051 --5 & 2.95 \\
  16 & 4.475 --2 & 1.00 & 3.599 --4 & 2.01 & 8.885 --6 & 2.99 \\
  32 & 2.238 --2 & 1.00 & 8.984 --5 & 2.00 & 1.113 --6 & 3.00 \\
  64 & 1.120 --2 & 1.00 & 2.245 --5 & 2.00 & 1.392 --7 & 3.00 \\
 128 & 5.599 --3 & 1.00 & 5.613 --6 & 2.00 & 1.740 --8 & 3.00 \\
 256 & 2.800 --3 & 1.00 & 1.403 --6 & 2.00 & 2.175 --9 & 3.00 \\
 512 & 1.400 --3 & 1.00 & 3.508 --7 & 2.00 & 2.719 --10 & 3.00 \\
1024 & 7.001 --4 & 1.00 & 8.769 --8 & 2.00 & 3.398 --11 & 3.00 \\
2048 & 3.501 --4 & 1.00 & 2.192 --8 & 2.00 & 4.236 --12 & 3.00 \\
\hline
\end{tabular}
\caption{Error for the regular function $u(t)=\sin \big( \frac{\pi}{4}t \big)$,
$t \in (0,T)$, with $T=2$.}\label{Table Fehler regular}
\end{table}

As a second example, we consider $u(t) = t^{2/3}$ for $t \in (0,T)$ and $T=2$ with
a singular behavior at $t=0$, i.e., $u \in H^{7/6-\varepsilon}_{0,}(0,T)$ for any $\varepsilon \in (0,1)$. Here, we observe
a convergence order of $\frac{2}{3}$ for $\nu \in \{0,1,2\}$, which for $\nu=1$, is between the expected order
$\frac{1}{6}$ and the optimal approximation order $\frac{7}{6}$, see
Table~\ref{Table Fehler singular 0}.

\begin{table}[ht]
\centering
\begin{tabular}{rcccccc}
\hline
 & \multicolumn{2}{c}{$\nu=0$} 
 & \multicolumn{2}{c}{$\nu=1$} 
 & \multicolumn{2}{c}{$\nu=2$} \\ \hline
$n$ & $\| u - u_h \|_{L^2(0,T)}$ & eoc
    & $\| u - u_h \|_{L^2(0,T)}$ & eoc
    & $\| u - u_h \|_{L^2(0,T)}$ & eoc \\ \hline
   2 & 4.796 --1 &      & 1.597 --1 &      & 6.268 --2 &      \\
   4 & 2.921 --1 & 0.72 & 9.222 --2 & 0.79 & 3.693 --2 & 0.76 \\
   8 & 1.783 --1 & 0.71 & 5.506 --2 & 0.74 & 2.237 --2 & 0.72 \\
  16 & 1.093 --1 & 0.71 & 3.369 --2 & 0.71 & 1.381 --2 & 0.70 \\
  32 & 6.741 --2 & 0.70 & 2.090 --2 & 0.69 & 8.604 --3 & 0.68 \\
  64 & 4.181 --2 & 0.69 & 1.306 --2 & 0.67 & 5.391 --3 & 0.67 \\
 128 & 2.605 --2 & 0.68 & 8.195 --3 & 0.67 & 3.387 --3 & 0.67 \\
 256 & 1.628 --2 & 0.68 & 5.152 --3 & 0.67 & 2.131 --3 & 0.67 \\
 512 & 1.021 --2 & 0.67 & 3.243 --3 & 0.67 & 1.341 --3 & 0.67 \\
1024 & 6.408 --3 & 0.67 & 2.042 --3 & 0.67 & 8.447 --4 & 0.67 \\
2048 & 4.028 --3 & 0.67 & 1.286 --3 & 0.67 & 5.320 --4 & 0.67 \\
\hline
\end{tabular}
\caption{Error for $u(t)=t^{2/3}$,
$t \in (0,T)$, with $T=2$.}\label{Table Fehler singular 0}
\end{table}

As a last example, we consider $u(t)=t(T-t)^{2/3}$, $t \in (0,T)$ with a
singular behavior at $t=T=2$. Again, we have 
$u \in H^{7/6-\varepsilon}_{0,}(0,T)$ but we observe the optimal order
of convergence as expected from the approximation property, see
Table \ref{Table Fehler singular T}.

\begin{table}[ht]
\centering
\begin{tabular}{rcccccc}
\hline
 & \multicolumn{2}{c}{$\nu=0$} 
 & \multicolumn{2}{c}{$\nu=1$} 
 & \multicolumn{2}{c}{$\nu=2$} \\ \hline
$n$ & $\| u - u_h \|_{L^2(0,T)}$ & eoc
    & $\| u - u_h \|_{L^2(0,T)}$ & eoc
    & $\| u - u_h \|_{L^2(0,T)}$ & eoc \\ \hline
   2 & 8.327 --1 &      & 1.506 --1 &      & 3.061 --2 &      \\
   4 & 3.980 --1 & 1.07 & 5.054 --2 & 1.58 & 1.221 --2 & 1.33 \\
   8 & 1.968 --1 & 1.02 & 1.823 --2 & 1.47 & 5.210 --3 & 1.23 \\
  16 & 9.852 --2 & 1.00 & 7.180 --3 & 1.34 & 2.271 --3 & 1.20 \\
  32 & 4.951 --2 & 0.99 & 3.000 --3 & 1.26 & 1.001 --3 & 1.18 \\
  64 & 2.490 --2 & 0.99 & 1.295 --3 & 1.21 & 4.433 --4 & 1.17 \\
 128 & 1.252 --2 & 0.99 & 5.677 --4 & 1.19 & 1.969 --4 & 1.17 \\
 256 & 6.287 --3 & 0.99 & 2.510 --4 & 1.18 & 8.760 --5 & 1.17 \\
 512 & 3.156 --3 & 0.99 & 1.114 --4 & 1.17 & 3.899 --5 & 1.17 \\
1024 & 1.583 --3 & 1.00 & 4.952 --5 & 1.17 & 1.736 --5 & 1.17 \\
2048 & 7.936 --4 & 1.00 & 2.204 --5 & 1.17 & 7.733 --6 & 1.17 \\
\hline
\end{tabular}
\caption{Error for $u(t)=t (T-t)^{2/3}$,
$t \in (0,T)$, with $T=2$.}\label{Table Fehler singular T}
\end{table}

\begin{remark}
When introducing the test space $W_h = \{ \int_0^{\cdot} v_h(s) ds \in H^1_{0,}(0,T) : v_h \in V_h\}$ with $V_h=S_h^\nu(0,T)$ given in \eqref{Vh} for $\nu \in \{0,1,2\}$, and when considering
the Galerkin--Petrov variational formulation to find $u_h \in V_h$ such that
\begin{equation*}
\langle u_h , {\mathcal{H}}_T w_h \rangle_{L^2(0,T)} =
\langle u , {\mathcal{H}}_T w_h \rangle_{L^2(0,T)} \quad 
\mbox{for all} \; w_h \in W_h ,
\end{equation*}
we prove, as in Remark~\ref{Bem:infsupH12}, the discrete inf-sup stability
condition
\[
\| u_h \|_{[H^{1/2}_{,0}(0,T)]'} \leq
\sup\limits_{0 \neq w_h \in W_h}
\frac{\langle u_h , {\mathcal{H}}_T w_h \rangle_{L^2(0,T)}}
{\| w_h \|_{H^{1/2}_{0,}(0,T)}} \quad
\mbox{for all} \; u_h \in V_h,
\]
which also ensures optimal error estimates for the approximate solution $u_h$.
However, and due to the applications in mind, our main interest is in the
numerical stability analysis of the Galerkin--Bubnov formulation
\eqref{Hilbert L2 Projektion}.
\end{remark}

\section{Discrete inf-sup stability condition in $S_h^0(0,T)$}
\label{Section:inf-sup}
For $n$ finite elements, we start by considering a given piecewise constant function
$u_h \in S_h^0(0,T)$, i.e.,
\[
  u_h(t) = \sum\limits_{i=1}^n u_i \psi_i^0(t), \quad t \in [0,T],
\]
with the norm
\[
  \| u_h \|^2_{L^2(0,T)} = \int_0^T [u_h(t)]^2 \, dt \, = \,
  h \, \sum\limits_{i=1}^n u_i^2 \, .
\]
In addition, we consider its Fourier series, see \eqref{Fourier series cos},
\[
  u_h(t) =
  \sum\limits_{k=0}^\infty \overline{u}_k
  \cos \left( \Big( \frac{\pi}{2} + k \pi \Big) \frac{t}{T} \right),
  \quad
  \overline{u}_k \, = \, \frac{2}{T} \int_0^T u_h(t) \,
  \cos \left( \Big( \frac{\pi}{2} + k \pi \Big) \frac{t}{T} \right) \, dt ,
\]
with the norm representation \eqref{Parseval},
\[
  \| u_h \|^2_{L^2(0,T)} = \frac{T}{2} \, \sum\limits_{k=0}^\infty
  \overline{u}_k^2 \, .
\]
It turns out that for $u_h \in S_h^0(0,T)$, it is sufficient to use a finite
sum of the Fourier coefficients $\overline{u}_k^2$ to define an
equivalent norm. Before we state this result, we need an
auxiliary lemma, for which we first define
\begin{equation}\label{Def xk}
  x_k := \left( \frac{\pi}{2} + k\pi \right) \frac{1}{2n} \quad
  \mbox{for} \; k \in {\mathbb{N}}_0 .
\end{equation}

\begin{lemma}
  For all $n\in\N$, and for all $k \in \N_0$, the equality
  \begin{equation}\label{sum cos2}
    \sum_{i=1}^n \cos^2 ((2i-1)x_k) \, = \, 
    \sum_{i=1}^n \cos^2 \left( \Big(
      \frac{\pi}{2}+k\pi\Big) \frac{2i-1}{2n} \right) \, = \, \frac{n}{2}
  \end{equation}
  holds true.
\end{lemma}
\proof{When using $\cos(2x) = 2 \cos^2 x - 1$, we write
  \[
    \sum_{i=1}^n \cos^2 \left( \Big(
      \frac{\pi}{2}+k\pi\Big) \frac{2i-1}{2n} \right)
    \, = \, \frac{1}{2} \sum_{i=1}^n \left[ 1 + \cos \left( \Big(
        \frac{\pi}{2}+k\pi\Big) \frac{2i-1}{n} \right) \right],
  \]
  and the assertion follows when using
  \cite[Equation~1.342,~4.]{GradshteynRyzhik2015}.}

\begin{lemma}
Let $M=n^2$. Then, the norm equivalence inequalities
\begin{equation}\label{norm equivalence n2}
  \frac{T}{2} \, \sum\limits_{k=0}^M \overline{u}_k^2 \leq
  \| u_h \|_{L^2(0,T)}^2 \leq 2 \, \left(
\frac{T}{2} \, \sum\limits_{k=0}^M \overline{u}_k^2
\right)
\end{equation}
hold true.
\end{lemma}
\proof{
While the lower estimate is trivial, to prove the upper
estimate, we first compute the Fourier coefficients
\begin{equation}
  \overline{u}_k
  \, = \, \frac{2}{T} \int_0^T u_h(t) \,
        \cos \left( \Big( \frac{\pi}{2} + k \pi \Big) \frac{t}{T} \right)
        \, dt 
  \, = \,
  \frac{2}{n} \, \frac{\sin x_k}{x_k} \, \sum\limits_{i=1}^n u_i
  \cos ((2i-1)x_k).
  \label{Def uk}
\end{equation} 
Using H\"older's inequality and \eqref{sum cos2}, we estimate
\begin{align*}
  \overline{u}_k^2
  &= \frac{4}{n^2} \, \frac{\sin^2 x_k}{x_k^2} \, 
        \left[ \sum\limits_{i=1}^n u_i
        \cos \left( \Big( \frac{\pi}{2} + k \pi \Big)
        \frac{2i-1}{2n} \right) \right]^2 \\
  &\leq \frac{4}{n^2} \, \frac{\sin^2 x_k}{x_k^2} \sum\limits_{i=1}^n u_i^2 \;
           \sum\limits_{i=1}^n 
           \cos^2 \left( \Big( \frac{\pi}{2} + k \pi \Big)
           \frac{2i-1}{2n} \right) \\
  &=\frac{2}{n} \, \frac{\sin^2 x_k}{x_k^2} \sum\limits_{i=1}^n u_i^2 \, = \,
           \frac{2}{T} \, \frac{\sin^2 x_k}{x_k^2} \,
           \| u_h \|^2_{L^2(0,T)} .
\end{align*}
Hence, we write
\[
  \| u_h \|^2_{L^2(0,T)}
  \, = \, \frac{T}{2} \sum\limits_{k=0}^\infty \overline{u}_k^2 
  \, \leq \, \frac{T}{2} \sum\limits_{k=0}^M \overline{u}_k^2 +
           \sum\limits_{k=M+1}^\infty \frac{\sin^2 x_k}{x_k^2} \, 
           \| u_h \|^2_{L^2(0,T)},
\]
and we estimate
\[
  \sum\limits_{k=M+1}^\infty \frac{\sin^2 x_k}{x_k^2} =
  \sum\limits_{k=M+1}^\infty
  \frac{\sin^2 \left( \Big( \frac{\pi}{2} + k \pi \Big) \frac{1}{2n}
  \right)}{[(\frac{\pi}{2}+k\pi)\frac{1}{2n}]^2}
  \leq \frac{4n^2}{\pi^2} \sum\limits_{k=M+1}^\infty
  \frac{1}{k^2} \, \leq \,
  \frac{n^2}{2} \, \sum\limits_{k=M+1}^\infty \frac{1}{k^2} \, .
\]         
With
\[
  \sum\limits_{k=M+1}^\infty \frac{1}{k^2} \leq
  \int_M^\infty \frac{1}{x^2} \, dx =
  \left[ - \frac{1}{x} \right]_{x=M}^\infty = \frac{1}{M},
\]
we further conclude
\[
  \| u_h \|^2_{L^2(0,T)} \leq
  \frac{T}{2} \sum\limits_{k=0}^M \overline{u}_k^2 +
           \frac{n^2}{2} \, \frac{1}{M} \, \| u_h \|_{L^2(0,T)}^2,
\]
and due to $M=n^2$ we obtain the assertion.}

\noindent
From the norm equivalence inequalities \eqref{norm equivalence n2}, we
observe that $u_h \equiv 0$ if $\overline{u}_k=0$ for all $k=0,\ldots,n^2$.
Thus, all coefficients $\overline{u}_k$ for $k > n^2$ have to be linear
dependent on the coefficients $\overline{u}_k$ for $k=0,\ldots,n^2$.
In the following lemma, we state this relation in more detail.
In particular, we prove that the Fourier coefficients 
$\overline{u}_k$ for $k=0,\ldots,n-1$ are sufficient to describe $u_h$.

\begin{lemma}
  The Fourier coefficients as given in \eqref{Def uk} satisfy the
  recurrence relations
  \begin{equation}\label{uk+2mn}
    \overline{u}_{k+2\mu n} = \frac{ 2 k + 1} {2k +1 + 4\mu n} \,
    \overline{u}_k \quad \mbox{for } \; k = 0,1,2,\ldots,2n-1, \;
    \mu \in \N,
  \end{equation}
  and
  \begin{equation}\label{uk+n}
    \overline{u}_{2n-1-k} = - \frac{2k+1} {4n-1-2k} \,
    \overline{u}_k \quad
    \mbox{for } \; k=0,1,2\ldots,n-1.
  \end{equation}
\end{lemma}
  \proof{The assertion follows from direct computations, we skip the details.}

\noindent
With this, we are in a position to rewrite the upper norm equivalence
inequality in \eqref{norm equivalence n2} in a more appropriate way.

\begin{corollary}
  For $u_h \in V_h = S_h^0(0,T)$, the estimate
\begin{equation*}
    \| u_h \|_{L^2(0,T)}^2 \leq \frac{T}{2} \,
    \frac{\pi^2}{3} \, \sum\limits_{k=0}^{n-1} \overline{u}_k^2
\end{equation*}
holds true, where the Fourier coefficients $\overline{u}_k$ are given as in
\eqref{Def uk}. 
\end{corollary}

\proof{We define
\begin{equation*}
  \gamma(k,n) := \sum\limits_{\mu=0}^\infty
  \frac{ (2 k + 1)^2} {(2k +1 + 4\mu n)^2} ,
\end{equation*}
satisfying, for $k=0,\ldots,2n-1$,
  \[
    1 \leq \gamma(k,n) =
    \sum\limits_{\mu=0}^\infty \frac{(2k+1)^2}{(2k+1+4\mu n)^2} = 
    \sum\limits_{\mu=0}^\infty \frac{(\frac{2k+1}{4n})^2}{(\frac{2k+1}{4n}+\mu)^2}\leq
    \sum\limits_{\mu=0}^\infty \frac{1}{(1+\mu)^2} = \frac{\pi^2}{6} , 
  \]
  where we use that the function $(0,1) \ni y \mapsto \frac{y^2}{(y+\mu)^2} \in \R$ is non-decreasing.
  Equation~\eqref{uk+2mn} gives
  \[
    \sum\limits_{k=0}^\infty \overline{u}_k^2
    = \sum\limits_{k=0}^{2n-1} \sum\limits_{\mu=0}^\infty
    \overline{u}^2_{k+2\mu n}
    = \sum\limits_{k=0}^{2n-1} \overline{u}_k^2
    \sum\limits_{\mu=0}^\infty \frac{ (2 k + 1)^2} {(2k +1 + 4\mu n)^2}
    = \sum\limits_{k=0}^{2n-1} \gamma(k,n) \, \overline{u}_k^2 .
  \]
  Next, we employ \eqref{uk+n} and the transformation
  $k=2n-1-j$ for $j=0,\ldots,n-1$ to conclude that
  \begin{align*}
    \sum\limits_{k=0}^{2n-1} \gamma(k,n) \, \overline{u}_k^2
    &=\sum\limits_{k=0}^{n-1} \gamma(k,n) \, \overline{u}_k^2 +
          \sum\limits_{k=n}^{2n-1} \gamma(k,n) \, \overline{u}_k^2  \\
    &=\sum\limits_{k=0}^{n-1} \gamma(k,n) \, \overline{u}_k^2 +
          \sum\limits_{j=0}^{n-1} \gamma(2n-1-j,n) \,
          \overline{u}_{2n-1-j}^2 \\
    &=\sum\limits_{k=0}^{n-1} \left[
          \gamma(k,n) +
          \gamma(2n-1-k,n) \,
          \frac{(2k+1)^2}{(4n-1-2k)^2} \right] \, \overline{u}_k^2 \\
    &\leq \sum\limits_{k=0}^{n-1} \Big[
          \gamma(k,n) +
          \gamma(2n-1-k,n) \,
          \Big] \, \overline{u}_k^2 \leq \frac{\pi^2}{3} 
          \sum\limits_{k=0}^{n-1} \overline{u}_k^2 .
  \end{align*}
  When using this within the norm representation
  \eqref{Parseval}, this gives the assertion.}

\begin{lemma}
  Let $n \in \N$ be given. For all $k,\ell=0,\ldots,n-1$, the discrete orthogonality
\begin{equation}\label{discrete orthogonality}
\sum\limits_{i=1}^n \sin ((2i-1)x_k) \sin ((2i-1) x_\ell) =
\frac{n}{2} \, \delta_{k\ell}
\end{equation}
holds true.
\end{lemma}

\proof{For $k = \ell$, the assertion is a simple consequence of
  \eqref{sum cos2}, due to $\sin^2 = 1 - \cos^2x$. For the remaining
  case $\ell \neq k$, we have, using
  $\sin x \sin y = \frac{1}{2}[\cos (x-y) - \cos (x+y)]$,
  \begin{multline*}
    \sum\limits_{i=1}^n \sin ((2i-1)x_k) \sin ((2i-1) x_\ell)
    = \sum\limits_{i=1}^n \sin \left( \Big(
       \frac{\pi}{2} + k \pi \Big) \frac{2i-1}{2n} \right)
       \sin \left( \Big(
       \frac{\pi}{2} + \ell \pi \Big) \frac{2i-1}{2n} \right) \\
    = \, \frac{1}{2} \sum\limits_{i=1}^n \left[ \cos \left( (k - \ell) \pi
       \frac{2i-1}{2n} \right) -
       \cos \left( ( k + \ell + 1) \pi  \frac{2i-1}{2n} \right) \right] ,
  \end{multline*}
  and it is sufficient to prove
  \[
    \sum\limits_{i=1}^n \cos \left( \nu \pi
       \frac{2i-1}{2n} \right) = 0
  \]
  for all $\nu=1,\ldots,2n-1$.
  Using \cite[Equation~1.341,~3.]{GradshteynRyzhik2015} yields
  \begin{multline*}
    \sum_{i=1}^n \cos \left( \nu \pi \frac{2i-1}{2n} \right) =
    \sum_{j=0}^{n-1} \cos \left( \nu \pi \frac{2j+1}{2n} \right) \\
    = \underbrace{\cos \left( \nu \pi \frac{1}{2n} +
        \nu \pi \frac{n-1}{2n} \right)}_{=0}
    \sin \left( \frac{\nu \pi }{2} \right)
    \frac{1}{\sin  \left( \nu \pi \frac{1}{2n} \right)} = 0
  \end{multline*}
  and thus, the assertion.}

\noindent
With this, we are in the position to state the main result of this section:

\begin{theorem}\label{Theorem inf-sup}
  For an arbitrary but fixed $u_h \in S_h^0(0,T)$, we define
  $w_h := Q_h {\mathcal{H}}_T^{-1}u_h$ as the unique solution
  of the variational formulation
  \begin{equation}\label{Definition wh}
    \langle w_h , v_h \rangle_{L^2(0,T)} =
    \langle {\mathcal{H}}_T^{-1} u_h , v_h \rangle_{L^2(0,T)}
    \quad \mbox{for all} \; v_h \in S_h^0(0,T) .
  \end{equation}
  Then, the inf-sup condition
  \begin{equation}\label{inf sup uh wh}
    c_S(u_h) \, 
    \| u_h \|_{L^2(0,T)} \leq
    \frac{\langle u_h , {\mathcal{H}}_T w_h \rangle_{L^2(0,T)}}
    {\| w_h \|_{L^2(0,T)}} 
  \end{equation}
  holds true with
  \begin{equation}\label{csuh}
    c_S(u_h) :=
    \frac{2\sqrt{3}}{\pi^2} \,
          \frac{16n^2-8n(2M+1)}{(4n-(2M+1))^2},
  \end{equation}
  where
  \begin{equation}\label{Assumption uh}
    M := \mbox{\rm arg min} \left \{ m \in {\mathbb{N}} :
    \| u_h \|^2_{L^2(0,T)} \leq \frac{T}{2} \, \frac{\pi^2}{3} \,
    \sum\limits_{k=0}^m \overline{u}_k^2 \right \} \leq n-1 
\end{equation}
with the Fourier coefficients $\overline{u}_k$
as given in \eqref{Def uk}. 
\end{theorem}

\proof{For a given $u_h \in S_h^0(0,T)$, we first define
\[
  w(t) := {\mathcal{H}}_T^{-1} u_h(t) =
  \sum\limits_{k=0}^\infty \overline{u}_k
  \sin \left( \Big( \frac{\pi}{2} + k \pi \Big) \frac{t}{T} \right),
\]
where the Fourier coefficients $\overline{u}_k$ are as in
\eqref{Def uk}. For this, we compute the piecewise constant
$L^2$ projection
\[
  w_h(t) = Q_h {\mathcal{H}}_T^{-1} u_h(t) = \sum\limits_{j=1}^n
  w_j \psi_j^0(t)
\]
as the unique solution of \eqref{Definition wh}, from which we conclude
\begin{equation}\label{Berechnung wj}
  w_j = \sum\limits_{k=0}^\infty \overline{u}_k \, \frac{\sin x_k}{x_k} \,
  \sin ((2j-1)x_k),
\end{equation}
and using \eqref{uk+2mn}, we write
\begin{equation}\label{wj+2mn}
  w_j = \sum\limits_{k=0}^{2n-1} \overline{u}_k \,
  \gamma(k,n) \, \frac{\sin x_k}{x_k} \,
  \sin ((2j-1)x_k) .
\end{equation}
Next, we use \eqref{uk+n} to conclude
\begin{equation}\label{wj+n}
  w_j =
\sum\limits_{k=0}^{n-1} \overline{u}_k \Big[
              \gamma(k,n) 
              - \gamma(2n-1-k,n)
              \frac{(2k+1)^2}{(4n-1-2k)^2}
              \Big] \, \frac{\sin x_k}{x_k} \,
              \sin ((2j-1)x_k) . 
\end{equation}
  Using the discrete orthogonality
  \eqref{discrete orthogonality}, it remains to compute
  \begin{align*}
    \| w_h& \|^2_{L^2(0,T)}
    = h \sum\limits_{j=1}^n w_j^2 \\
    =& h \sum\limits_{j=1}^n \left \{
          \sum\limits_{k=0}^{n-1} \overline{u}_k \Big[
              \gamma(k,n) 
              - \gamma(2n-1-k,n)
              \frac{(2k+1)^2}{(4n-1-2k)^2}
              \Big] \frac{\sin x_k}{x_k}
              \sin ((2j-1)x_k) 
       \right \}^2 \\
    =& \frac{T}{2} \,
       \sum\limits_{k=0}^{n-1} \overline{u}_k^2 \, \Big[
       \gamma(k,n) - \gamma(2n-1-k,n) \frac{(2k+1)^2}{(4n-1-2k)^2}
       \Big]^2 \; \frac{\sin^2 x_k}{x_k^2}.
  \end{align*}
  With $2k+1 < 4n-(2k+1)$ for $k=0,\ldots,n-1$, we calculate that
  \begin{align*}
       \gamma(k,n) - \gamma(2n-1-k&,n) \frac{(2k+1)^2}{(4n-1-2k)^2} \\
    &= \sum\limits_{\mu=0}^\infty \frac{(2k+1)^2}{(2k+1+4\mu n)^2}
       - \frac{(2k+1)^2}{(4n-1-2k)^2} \sum\limits_{\mu=0}^\infty
       \frac{(4n-1-2k)^2}{(4n-1-2k+4\mu n)^2} \\
    &= (2k+1)^2 
     \sum\limits_{\mu=0}^\infty \left[ \frac{1}{(2k+1+4\mu n)^2}
     - 
     \frac{1}{(4n-1-2k+4\mu n)^2} \right] \\
    &> (2k+1)^2 
     \left[ \frac{1}{(2k+1)^2}
     - 
     \frac{1}{(4n-1-2k)^2} \right] \\
    &= 1 - \frac{(2k+1)^2}{(4n-1-2k)^2} .
  \end{align*}
  For $k=0,\ldots,n-1$, we have $x_k \in (0,\frac{\pi}{2})$, and therefore,
  \[
    \frac{\sin^2 x_k}{x_k^2} \geq
    \frac{4}{\pi^2} \quad \mbox{for all} \; k=0,\ldots,n-1.
  \]
  Let $M \leq n-1$ such that \eqref{Assumption uh} is satisfied. Then, we write
  \begin{align*}
    \| w_h \|^2_{L^2(0,T)}
    &> \frac{T}{2} \, \frac{4}{\pi^2} \,
          \sum\limits_{k=0}^{n-1} \overline{u}_k^2 \left[
          1 - \frac{(2k+1)^2}{(4n-1-2k)^2} \right]^2 \\
    &\geq \frac{T}{2} \, \frac{4}{\pi^2} \,
          \sum\limits_{k=0}^M \overline{u}_k^2 \left[
          1 - \frac{(2k+1)^2}{(4n-1-2k)^2} \right]^2 \\
    &\geq \frac{T}{2} \, \frac{4}{\pi^2} \,
             \left[
          1 - \frac{(2M+1)^2}{(4n-1-2M)^2} \right]^2
          \sum\limits_{k=0}^M \overline{u}_k^2 \\
    &\geq \frac{4}{\pi^2} \, \frac{3}{\pi^2}
             \left[
          1 - \frac{(2M+1)^2}{(4n-1-2M)^2} \right]^2 \, \| u_h \|^2_{L^2(0,T)} \\
    &= \frac{12}{\pi^4} \, \left[
          \frac{16n^2-8n(2M+1)}{(4n-(2M+1))^2} \right]^2
          \, \| u_h \|^2_{L^2(0,T)} ,
  \end{align*}
  i.e., we conclude that
  \[
    \| w_h \|_{L^2(0,T)} \, \geq \,
    \frac{2\sqrt{3}}{\pi^2} \,
          \frac{16n^2-8n(2M+1)}{(4n-(2M+1))^2}
    \, \| u_h \|_{L^2(0,T)}.
  \]
  Due to
  \[
    \langle u_h , {\mathcal{H}}_T w_h \rangle_{L^2(0,T)}
    \, = \, \langle {\mathcal{H}}_T^{-1} u_h , w_h \rangle_{L^2(0,T)} \, = \,
          \langle w_h , w_h \rangle_{L^2(0,T)} 
    \, = \, \| w_h \|_{L^2(0,T)}^2 ,
  \]
  this finally implies the desired estimate.}

\begin{corollary}
  The estimate \eqref{inf sup uh wh} implies the inf-sup stability
  condition
  \[
    c_S \, 
    \| u_h \|_{L^2(0,T)} \leq \sup\limits_{0 \neq v_h \in S_h^0(0,T)}
    \frac{\langle u_h , {\mathcal{H}}_T v_h \rangle_{L^2(0,T)}}
    {\| v_h \|_{L^2(0,T)}}
  \]
  for all $u_h \in S_h^0(0,T)$, where we have to consider $M=n-1$,
  i.e.,
  \[
    c_S := \min\limits_{u_h\in S_h^0(0,T)} c_S(u_h) =
    \frac{1}{T} \,
    \frac{2\sqrt{3}}{\pi^2} \,
          \frac{8}{(2+\frac{1}{n})^2} \, h \, .
  \]
  In particular for $T=2$, we have
  \[
    \lim\limits_{n \to \infty}
    \frac{1}{2} \,
    \frac{2\sqrt{3}}{\pi^2} \,
    \frac{8}{(2+\frac{1}{n})^2} \, = \,
    \frac{2\sqrt{3}}{\pi^2} \, \simeq \, 0.351 .
  \]
  Note that the value obtained from the numerical
  experiments as given in Table \ref{Table inf sup constant L2} was
  $0.426$.
\end{corollary}

\begin{remark}
  In the case of a constant $M$ independent of $n$, we obtain
  \[
    c_S(u_h) \simeq \frac{2\sqrt{3}}{\pi^2}
  \]
  independent of the mesh size $h$. In the case $M=n-\sqrt{n}$, this
  would result in
  \[
    c_S(u_h) = \frac{2\sqrt{3}}{\pi^2} \, \frac{16n^2-8n(2M+1)}{(4n-(2M+1))^2}
    =
    \frac{2\sqrt{3}}{\pi^2} \,
    \frac{16n\sqrt{n}-1}{(2n+2\sqrt{n}-1)^2}
    \simeq \frac{8\sqrt{3}}{\pi^2} \, \frac{1}{\sqrt{n}} =
    \frac{8\sqrt{3}}{\pi^2 \sqrt{T}} \, h^{1/2} \, .
  \]
\end{remark}

\section{Error estimates for the projection in $S_h^0(0,T)$}
\label{Section:Error}
This section aims to present a numerical analysis in order
to confirm the numerical results as presented in
Section \ref{Section:Projection}. In addition to the solution
$u_h \in V_h = S_h^0(0,T)$ of the variational formulation
\eqref{Hilbert L2 Projektion}, we introduce the $L^2$ projection
$Q_hu \in S_h^0(0,T)$ as unique solution of the variational formulation
\begin{equation}\label{Def L2 Projektion}
  \langle Q_h u , v_h \rangle_{L^2(0,T)} =
  \langle u , v_h \rangle_{L^2(0,T)} \quad \mbox{for all} \;
  v_h \in S_h^0(0,T) .
\end{equation}
When we assume $u \in H^s(0,T)$ for some $s \in [0,1]$, the
error estimate
\begin{equation*}
\| u - Q_h u \|_{L^2(0,T)} \leq c \, h^s \, \|u\|_{H^s(0,T)}
\end{equation*}
holds true, see, e.g., \cite{Steinbach:2008}.
Using the triangle inequality, we therefore have
\begin{equation}\label{Error u uh triangle}
  \| u - u_h \|_{L^2(0,T)} \leq c \, h^s \, \|u\|_{H^s(0,T)} +
  \| u_h - Q_h u \|_{L^2(0,T)},
\end{equation}
and it remains to estimate the second term.

\begin{lemma}
  Let $u_h \in S_h^0(0,T)$ be the unique solution of the variational
  formulation \eqref{Hilbert L2 Projektion}, and let $Q_hu \in S_h^0(0,T)$
  be the $L^2$ projection of $u \in L^2(0,T)$ as defined in 
  \eqref{Def L2 Projektion}. Then, the estimate
  \begin{equation}\label{estimate uh-Qhu}
    \| u_h - Q_h u \|_{L^2(0,T)} \leq \frac{1}{c_S(u_h-Q_hu)} \,
    \| Q_h {\mathcal{H}}_T^{-1} (u-Q_hu) \|_{L^2(0,T)}
  \end{equation}
  holds true,
  where $c_S(u_h-Q_hu)$ is the stability constant as defined in
  \eqref{csuh}.
\end{lemma}

\proof{As in Theorem \ref{Theorem inf-sup},
we define $w={\mathcal{H}}_T^{-1}(u_h-Q_hu)$ as well as
$w_h=Q_h {\mathcal{H}}_T^{-1}(u_h-Q_hu)$. From the definition of
$w_h \in S_h^0(0,T)$ and relation~\eqref{HT-HTinverse}, we have
\begin{align*}
  \langle w_h , v_h \rangle_{L^2(0,T)}
  &= \langle {\mathcal{H}}_T^{-1}(u_h-Q_hu) , v_h \rangle_{L^2(0,T)}
        \, = \,
        \langle u_h - Q_h u , {\mathcal{H}}_T v_h \rangle_{L^2(0,T)} \\
  &= \langle u - Q_h u , {\mathcal{H}}_T v_h \rangle_{L^2(0,T)} \, = \,
        \langle {\mathcal{H}}_T^{-1}(u -Q_hu) , v_h \rangle_{L^2(0,T)}
\end{align*}
for all $v_h\in S_h^0(0,T)$,
i.e., $w_h =Q_h {\mathcal{H}}_T^{-1}(u-Q_hu)$. Hence, using \eqref{Hilbert L2 Projektion}, we write
\eqref{inf sup uh wh} as
\begin{align*}
  c_S(u_h-Q_hu) \, &\| u_h - Q_hu \|_{L^2(0,T)} \leq
     \frac{\langle u_h - Q_h u , {\mathcal{H}}_T w_h \rangle_{L^2(0,T)}}
     {\| w_h \|_{L^2(0,T)}} \\
  &= \frac{\langle u - Q_h u , {\mathcal{H}}_T w_h \rangle_{L^2(0,T)}}
     {\| w_h \|_{L^2(0,T)}} \, = \,
     \frac{\langle {\mathcal{H}}_T^{-1}(u - Q_h u) , w_h \rangle_{L^2(0,T)}}
     {\| w_h \|_{L^2(0,T)}} \\
  &=\frac{\langle w_h , w_h \rangle_{L^2(0,T)}}
     {\| w_h \|_{L^2(0,T)}} \, = \, \| w_h \|_{L^2(0,T)} \, = \,
     \| Q_h {\mathcal{H}}_T^{-1} (u-Q_hu) \|_{L^2(0,T)} ,
\end{align*}
i.e., the assertion follows.}

\begin{remark}
When using the stability of 
$Q_h \colon \, L^2(0,T) \to S_h^0(0,T) \subset L^2(0,T)$ and Parseval's theorem
for the inverse ${\mathcal{H}}_T^{-1}$ of the modified Hilbert transformation,
the estimates
\begin{equation}\label{Estimate cB H1}
  \| Q_h {\mathcal{H}}_T^{-1} (u-Q_hu) \|_{L^2(0,T)} \leq
  \| u -Q_h u \|_{L^2(0,T)} \leq c \, h^s \, \|u\|_{H^s(0,T)}
\end{equation}
hold true for $s \in [0,1]$. Combining these estimates with \eqref{Error u uh triangle} and
\eqref{estimate uh-Qhu} yields
\[
  \| u - u_h \|_{L^2(0,T)} \leq c \, \left( 1 + 
  \frac{1}{c_S(u_h-Q_hu)} \right) \, h^s \, \|u\|_{H^s(0,T)}
\]
for $s \in [0,1]$, which only results in optimal convergence rates when $c_S(u_h-Q_hu)=
{\mathcal{O}}(1)$, which we can not expect to hold in general.
\end{remark}

\noindent
Hence, we have to analyze $\| Q_h {\mathcal{H}}_T^{-1} (u-Q_hu) \|_{L^2(0,T)}$
in more detail in order to prove second-order convergence when assuming
$u \in H^2(0,T)$. We define
\[
U(t) := \int_0^t u(s) \, ds, \quad U(0)=0, \quad \partial_t U(t) = u(t).
\]
The coefficients $u_i$ of the piecewise constant $L^2$ projection $Q_hu$
are then given by
\[
  u_i = \frac{1}{h} \int_{t_{i-1}}^{t_i} u(s) \, ds =
  \frac{1}{h} \int_{t_{i-1}}^{t_i} \partial_t U(s) \, ds =
  \frac{1}{h} \Big[ U(t_i) - U(t_{i-1}) \Big] = (\partial_t I_hU)_{|(t_{i-1},t_i)}(t),
\]
where $I_hU \in S_h^1(0,T)$ is the piecewise linear interpolation of $U$ in 
$S_h^1(0,T)$, see \eqref{Vh}. Hence,
$Q_hu= \partial_t I_hU$ in $(t_{i-1},t_i)$ for all $i=1,\ldots,n$, and
in particular we have $u-Q_hu = \partial_t(U-I_hU)$.

\begin{lemma}
  For the interpolation error, the local
  representation
  \begin{equation*}
    U(t) - I_hU(t) = \int_{t_{i-1}}^{t_i} G(s,t) \, \partial_{tt} U(s) \, ds, \quad t \in (t_{i-1},t_i),
  \end{equation*}
  holds true with Green's function
  \[
    G(s,t) = \frac{1}{h} \left \{
      \begin{array}{ccl}
        (s-t_{i-1})(t-t_i), && \; s \in (t_{i-1},t), \\[1mm]
        (s-t_i)(t-t_{i-1}), && \; s \in (t,t_i) .
      \end{array}
      \right.
  \]
\end{lemma}

\proof{Although interpolation error estimates are well-known,
  for completeness, we give a proof for this particular result.
  Recall that for $t \in (t_{i-1},t_i)$ the piecewise linear interpolation
  is given as
  \[
    I_hU(t) = U(t_{i-1}) + \frac{1}{h} (t-t_{i-1})
    \Big[ U(t_i)-U(t_{i-1}) \Big].  
  \]
  On the other hand, for $t \in (t_{i-1},t_i)$, we calculate that
  \begin{align*}
    U(t) - U(t_{i-1})
    &= \int_{t_{i-1}}^t \partial_t U(s) \, ds \, = \,
          (s-t) \partial_t U(s) \Big|_{t_{i-1}}^t -
          \int_{t_{i-1}}^t (s-t) \, \partial_{tt} U(s) \, ds \\
    &= (t-t_{i-1}) \partial_t U(t_{i-1}) + \int_{t_{i-1}}^t (t-s) \, \partial_{tt} U(s) \, ds \, .
  \end{align*}
  In particular, for $t=t_i$, we have
  \[
    U(t_i) - U(t_{i-1}) \, = \,
    h \, \partial_t U(t_{i-1}) + \int_{t_{i-1}}^{t_i} (t_i-s) \, \partial_{tt} U(s) \, ds \, ,
  \]
  i.e.,
  \[
    \partial_t U(t_{i-1}) \, = \, \frac{1}{h} \Big[ U(t_i)-U(t_{i-1}) \Big] - \frac{1}{h}
    \int_{t_{i-1}}^{t_i} (t_i-s) \, \partial_{tt} U(s) \, ds .
  \]
  Hence, we obtain
  \begin{multline*}
    U(t)
    = U(t_{i-1}) + \frac{1}{h} (t-t_{i-1}) \Big[ U(t_i)-U(t_{i-1}) \Big] \\
       - \frac{1}{h}(t-t_{i-1}) \int_{t_{i-1}}^{t_i} (t_i-s) \, \partial_{tt} U(s) \, ds
       + \int_{t_{i-1}}^t (t-s) \, \partial_{tt} U(s) \, ds,
  \end{multline*}
  and
  \begin{align*}
    U(t) - I_hU(t)
    & =\int_{t_{i-1}}^t (t-s) \, \partial_{tt} U(s) \, ds -
          \frac{1}{h}(t-t_{i-1}) \int_{t_{i-1}}^{t_i} (t_i-s) \, \partial_{tt} U(s) \, ds \\
    & = \frac{1}{h} (t-t_i) \int_{t_{i-1}}^t (s-t_{i-1}) \, \partial_{tt} U(s) \, ds 
          +
          \frac{1}{h}(t-t_{i-1}) \int_t^{t_i} (s-t_i) \, \partial_{tt} U(s) \, ds .
  \end{align*}
This concludes the proof.}

\noindent
With $\partial_{tt} U(t)=\partial_t u(t)$ we therefore have, using integration by parts,
\begin{align*}
  U(t) - &I_hU(t)
  = \frac{1}{h} (t-t_i) \int_{t_{i-1}}^t (s-t_{i-1}) \, \partial_t u(s) \, ds 
        +
        \frac{1}{h}(t-t_{i-1}) \int_t^{t_i} (s-t_i) \, \partial_t u(s) \, ds \\
  =& \frac{1}{h}(t-t_i) \left[
     \left. \frac{1}{2} \Big(
     (s-t_{i-1})^2 - h(t-t_{i-1}) \Big) \partial_t u(s) \right|_{t_{i-1}}^t \right. \\
  & \left. - \frac{1}{2} \int_{t_{i-1}}^t
     \Big( (s-t_{i-1})^2 - h(t-t_{i-1}) \Big) \, \partial_{tt} u(s) \, ds \right] \\
  & + \frac{1}{h}(t-t_{i-1}) \left[
     \left. \frac{1}{2} (s-t_i)^2 \, \partial_t u(s) \right|_t^{t_i} -
     \frac{1}{2} \int_t^{t_i} (s-t_i)^2 \, \partial_{tt} u(s) \, ds \right] \\
  =& \frac{1}{2} (t-t_i)(t-t_{i-1}) \partial_t u(t_{i-1}) 
     - \frac{1}{2} \frac{1}{h}(t-t_i) \int_{t_{i-1}}^t
     \Big( (s-t_{i-1})^2 - h(t-t_{i-1}) \Big) \, \partial_{tt} u(s) \, ds \\
  & - \frac{1}{2} \frac{1}{h} (t-t_{i-1})
     \int_t^{t_i} (s-t_i)^2 \, \partial_{tt} u(s) \, ds \, .
\end{align*}
Due to $u(t) - Q_hu(t) = \partial_t (U(t)-I_hU(t)) = u^1(t)+u^2(t)$,
we define
\begin{equation}\label{Def u1}
  u^1(t) := \partial_t \left[
    \frac{1}{2} (t-t_i)(t-t_{i-1}) \partial_t u(t_{i-1}) \right] =
  \left( t - \frac{1}{2} \Big( t_i+t_{i-1} \Big) \right) \, \partial_t u(t_{i-1}),
\end{equation}
and
\begin{align}
  u^2(t) \nonumber
  &:= - \frac{1}{2} \frac{1}{h} \partial_t \left[
        (t-t_i) \int_{t_{i-1}}^t
        \Big( (s-t_{i-1})^2 - h(t-t_{i-1}) \Big) \, \partial_{tt} u(s) \, ds \right. \\
  & \qquad \qquad \qquad \qquad \qquad \qquad \qquad \qquad \left. 
     + (t-t_{i-1}) \int_t^{t_i} (s-t_i)^2 \, \partial_{tt} u(s) \, ds \right] \nonumber \\
  &=- \frac{1}{2} \frac{1}{h} \left[
      \int_{t_{i-1}}^t \Big( (s-t_{i-1})^2 - h (2t-t_{i-1}-t_i) \Big) \,
      \partial_{tt} u(s) \, ds + \int_t^{t_i} (s-t_i)^2 \, \partial_{tt} u(s) \, ds
      \right] \nonumber \\
  &= - \frac{1}{2} \frac{1}{h}
        \int_{t_{i-1}}^{t_i} \widetilde{G}(s,t) \, \partial_{tt} u(s) \, ds
        \label{Def u2}
\end{align}
for $t \in (0,T)$ with the function
\[
  \widetilde{G}(s,t) = \left \{
    \begin{array}{ccl}
      (s-t_{i-1})^2 - h(2t-t_{i-1}-t_i)
      & & \mbox{for} \; s \in (t_{i-1},t), \\[1mm]
      (s-t_i)^2 & & \mbox{for} \; s \in (t,t_i).
    \end{array} \right.
\]
With this splitting, we have
\begin{equation}\label{Splitting}
  \| Q_h {\mathcal{H}}_T^{-1} (u-Q_hu) \|_{L^2(0,T)} \leq
  \| Q_h {\mathcal{H}}_T^{-1} u^1 \|_{L^2(0,T)}
  +
  \| u^2 \|_{L^2(0,T)} 
\end{equation}
where in the second argument, we used the boundedness of $Q_h$ and
${\mathcal{H}}_T^{-1}$. We show that both terms are of
order $h^2$ when we assume $u \in H^2(0,T)$.

\begin{lemma}
  Assume $u \in H^2(0,T)$, and let $u^2$ be as defined in
  \eqref{Def u2}. Then,
  \begin{equation}\label{Norm u2}
    \| u^2 \|_{L^2(0,T)} \leq \frac{1}{3} \, h^2 \, \| \partial_{tt} u \|_{L^2(0,T)} .
  \end{equation}
\end{lemma}
\proof{From \eqref{Def u2}, we immediately have
  \[
    [u^2(t)]^2 = \frac{1}{4} \frac{1}{h^2} \left(
      \int_{t_{i-1}}^{t_i} \widetilde{G}(s,t) \, \partial_{tt} u(s) \, ds \right)^2 \leq
    \frac{1}{4} \frac{1}{h^2} \int_{t_{i-1}}^{t_i}
    \Big[ \widetilde{G}(s,t) \Big]^2 \, ds
    \int_{t_{i-1}}^{t_i} [\partial_{tt} u(s)]^2 \, ds
  \]
  for $t \in (t_{i-1},t_i)$, i.e.,
  \[
    \int_{t_{i-1}}^{t_i} [u^2(t)]^2 \, dt \leq
    \frac{1}{4} \frac{1}{h^2} \int_{t_{i-1}}^{t_i} \int_{t_{i-1}}^{t_i}
    \Big[ \widetilde{G}(s,t) \Big]^2 \, ds \, dt
    \int_{t_{i-1}}^{t_i} [\partial_{tt} u(s)]^2 \, ds.
  \]
  A direct computation gives
  \[
    \int_{t_{i-1}}^{t_i} \int_{t_{i-1}}^{t_i}
    \Big[ \widetilde{G}(s,t) \Big]^2 \, ds \, dt =
    \frac{3}{10} \, h^6 ,
  \]
  and hence,
  \[
    \int_{t_{i-1}}^{t_i} [u^2(t)]^2 \, dt \leq
    \frac{3}{40} \, h^4 \int_{t_{i-1}}^{t_i} [\partial_{tt} u(s)]^2 \, ds .
  \]
  Summing up for $i=1,\ldots,n$, using the simple estimate
  $\frac{3}{40} < \frac{4}{36}$, and taking the square root,
  this gives the assertion.}

\noindent
It remains to estimate $Q_h {\mathcal{H}}_T^{-1}u^1$.
For this, we first consider the Fourier series
\[
  u^1(t) = \sum\limits_{k=0}^\infty \overline{u}_k^1
  \cos \left( \Big( \frac{\pi}{2} + k \pi \Big) \frac{t}{T} \right)
\]
with the Fourier coefficients, recall \eqref{Def xk} for the definition of
$x_k$,
\begin{align}
  \overline{u}_k^1
  &= \frac{2}{T} \int_0^T u^1(t) \,
        \cos \left( \Big( \frac{\pi}{2} + k \pi \Big) \frac{t}{T} \right)
        \, dt  \nonumber \\ 
  &= \frac{2}{T} \sum\limits_{i=1}^n \partial_t u(t_{i-1}) \int_{t_{i-1}}^{t_i}
        \left( t - \frac{1}{2} \Big( t_i + t_{i-1} \Big) \right)
        \cos \left( \Big( \frac{\pi}{2} + k \pi \Big) \frac{t}{T} \right)
        \, dt \nonumber \\
  &= \frac{1}{x_k} \frac{h}{n} \left[
        \cos x_k - \frac{\sin x_k}{x_k} \right]
        \sum\limits_{i=1}^n \partial_t u(t_{i-1}) \sin ((2i-1)x_k) .
        \label{Def uk1}
\end{align}
Hence, we have
\[
  {\mathcal{H}}_T^{-1} u^1(t) = \sum\limits_{k=0}^\infty \overline{u}_k^1
  \sin \left( \Big( \frac{\pi}{2} + k \pi \Big) \frac{t}{T} \right),
\]
and we compute the coefficients $w_i^1$ of
$w_h^1 = Q_h {\mathcal{H}}_T^{-1} u^1$ for $i=1,\ldots,n$, resulting in
\begin{align*}
  w_i^1 &= \frac{1}{h} \int_{t_{i-1}}^{t_i} {\mathcal{H}}_T^{-1} u^1(t) \, dt
            \, = \, \frac{1}{h} \sum\limits_{k=0}^\infty
            \overline{u}_k^1 \int_{t_{i-1}}^{t_i}
            \sin \left( \Big( \frac{\pi}{2} + k \pi \Big) \frac{t}{T} \right)
            \, dt \\
      &=\sum\limits_{k=0}^\infty \overline{u}_k^1 \, \frac{\sin x_k}{x_k}
            \, \sin ((2i-1)x_k) \, .
\end{align*}
For $k=0,\ldots,2n-1$ and $\mu \in {\mathbb{N}}_0$, we write
the Fourier coefficients \eqref{Def uk1} as
\begin{align*}
  \overline{u}_{k+2\mu n}^1
  &= \frac{1}{x_{k+2\mu n}} \frac{h}{n}
        \left[ \cos x_{k+2\mu n} - \frac{\sin x_{k+2\mu n}}{x_{k+2\mu n}} \right]
        \sum\limits_{i=1}^n \partial_t u(t_{i-1}) \, \sin ((2i-1) x_{k+2\mu n}) \\
  &= \frac{1}{(x_k+\mu\pi)^2} \frac{h}{n}
        \Big[ (x_k+\mu \pi) \cos x_k - \sin x_k \Big]
        \sum\limits_{i=1}^n \partial_t u(t_{i-1}) \, \sin ((2i-1) x_k) \\
  &=\left( \frac{x_k}{x_k+\mu \pi} \right)^2
        \frac{(x_k+\mu \pi) \cos x_k - \sin x_k}{x_k \cos x_k-\sin x_k} \,
        \overline{u}_k^1 \, .
\end{align*}
Thus, we compute that
\begin{align*}
  w_i^1
  &= \sum\limits_{k=0}^\infty \, \overline{u}_k^1 \,
        \frac{\sin x_k}{x_k} \, \sin ((2i-1) x_k) \\
  &= \sum\limits_{k=0}^{2n-1} \sum\limits_{\mu=0}^\infty
        \overline{u}_{k+2\mu n}^1 \,
        \frac{\sin x_{k+2\mu n}}{x_{k+2\mu n}} \, \sin ((2i-1) x_{k+2\mu n}) \\
  &= \sum\limits_{k=0}^{2n-1} \sum\limits_{\mu=0}^\infty
        \left( \frac{x_k}{x_k+\mu \pi} \right)^3
        \frac{(x_k+\mu \pi) \cos x_k - \sin x_k}{x_k \cos x_k-\sin x_k} \,
         \frac{\sin x_k}{x_k} \, \overline{u}_k^1 \, \sin ((2i-1) x_k) \\
  &= \sum\limits_{k=0}^{2n-1} \beta_k \, 
        \frac{\sin x_k}{x_k} \, \overline{u}_k^1 \, \sin ((2i-1) x_k),
\end{align*}
with
\[
  \beta_k
  = \sum\limits_{\mu=0}^\infty \left( \frac{x_k}{x_k+\mu \pi} \right)^3
  \frac{(x_k+\mu \pi) \cos x_k - \sin x_k}{x_k \cos x_k-\sin x_k} .
\]
For $k=0,\ldots,n-1$, we further obtain
\begin{align*}
  \overline{u}_{2n-1-k}^1
  &= \frac{1}{x_{2n-1-k}} \frac{h}{n}
        \left[ \cos x_{2n-1-k} - \frac{\sin x_{2n-1-k}}{x_{2n-1-k}} \right]
       \sum\limits_{i=1}^n \partial_t u(t_{i-1})
        \sin ((2i-1) x_{2n-1-k}) \\
  &= - \frac{1}{(\pi-x_k)^2} 
          \Big[ (\pi-x_k) \cos x_k + \sin x_k \Big]
          \frac{h}{n}
       \sum\limits_{i=1}^n \partial_t u(t_{i-1})
          \sin ((2i-1) x_k) \\
  &= - \frac{x_k^2}{(\pi-x_k)^2} 
            \frac{(\pi-x_k) \cos x_k + \sin x_k}{x_k \cos x_k - \sin x_k}
            \overline{u}_k^1 \, .
\end{align*}
Hence, this yields
\begin{align} \nonumber
  w_i^1 &= \sum\limits_{k=0}^{2n-1} \beta_k \,
        \frac{\sin x_k}{x_k} \, \overline{u}_k^1 \, \sin ((2i-1) x_k) \\
  =& \sum\limits_{k=0}^{n-1} \beta_k \, \nonumber
        \frac{\sin x_k}{x_k} \, \overline{u}_k^1 \, \sin ((2i-1) x_k)
        + \sum\limits_{k=0}^{n-1} \beta_{2n-1-k} \,
        \frac{\sin x_{2n-1-k}}{x_{2n-1-k}} \,
        \overline{u}_{2n-1-k}^1 \, \sin ((2i-1) x_{2n-1-k}) \\
  =& \sum\limits_{k=0}^{n-1} \left[ \beta_k - \beta_{2n-1-k} \,
        \frac{x_k^3}{(\pi-x_k)^3} \nonumber
        \frac{(\pi-x_k) \cos x_k + \sin x_k}{x_k \cos x_k - \sin x_k}
        \right] \frac{\sin x_k}{x_k} \,
        \overline{u}_k^1 \, \sin ((2i-1) x_k) \\
  =& \sum\limits_{k=0}^{n-1} F(x_k) \, \overline{u}_k^1 \,
        \sin ((2i-1)x_k) \label{Darstellung wi1}
\end{align}
with the function
\begin{multline*}
  F(x_k) = \left[ \beta_k 
        - \beta_{2n-1-k} \,
        \frac{x_k^3}{(\pi-x_k)^3} 
        \frac{(\pi-x_k) \cos x_k + \sin x_k}{x_k \cos x_k - \sin x_k}
             \right] \frac{\sin x_k}{x_k} \\
   = \left[ 1 + \frac{x_k^3}{x_k \cos x_k - \sin x_k} 
        \sum\limits_{\mu=1}^\infty \left(
        \frac{(x_k+\mu \pi) \cos x_k - \sin x_k}{(x_k+\mu \pi)^3}
        -
        \frac{(\mu\pi-x_k) \cos x_k +
        \sin x_k}{(\mu\pi-x_k)^3}
        \right) \,
     \right] \\
     \cdot\frac{\sin x_k}{x_k} .          
\end{multline*}

\begin{lemma}
  For a given $u \in H^2(0,T)$, let $u^1$ as defined in
  \eqref{Def u1}, i.e.,
  \[
    u^1(t) = \left( t - \frac{1}{2} \Big( t_i + t_{i-1} \Big) \right)
    \partial_t u(t_{i-1}) \quad \mbox{for} \; t \in (t_{i-1},t_i), \;
    i=1,\ldots,n.
  \]
  Define $w_h^1 := Q_h {\mathcal{H}}_T^{-1}u^1 \in S_h^0(0,T)$. Then,
  the estimate
  \begin{equation}\label{Estimate wh1}
    \| w_h^1 \|^2_{L^2(0,T)} \leq
    \frac{\pi}{96} \, h^4 \, \| \partial_{tt} u \|^2_{L^2(0,T)} +
    \frac{\pi}{48} \, h^3 \, [\partial_t u(0)]^2
  \end{equation}
  holds true. In particular, for $u\in H^2(0,T)$ with $\partial_t u(0)=0$, this gives
  \begin{equation}\label{Estimate wh10}
    \| w_h^1 \|^2_{L^2(0,T)} \leq
    \frac{\pi}{96} \, h^4 \, \| \partial_{tt} u \|^2_{L^2(0,T)}  \, .
  \end{equation}
\end{lemma}
\proof{
  Since $w_h^1 \in S_h^0(0,T)$ is piecewise constant, we write,
  using \eqref{Darstellung wi1} and \eqref{Def uk1},
  \begin{align*}
    \| w_h^1 \|^2_{L^2(0,T)}
    &= \int_0^T [w_h^1(t)]^2 \, dt \, = \,
          h \sum\limits_{i=1}^n [w_i^1]^2 \, = \,
          h \sum\limits_{i=1}^n \left[ \sum\limits_{k=0}^{n-1}
          F(x_k) \, \overline{u}_k^1 \, \sin ((2i-1) x_k) \right]^2 \\
    &= h \sum\limits_{k=0}^{n-1} \sum\limits_{\ell=0}^{n-1}
          F(x_k) \, F(x_\ell) \, \overline{u}_k^1 \, \overline{u}_\ell^1
          \sum\limits_{i=1}^n \Big[ \sin ((2i-1)x_k \sin ((2i-1)x_\ell) \Big] \\
    &= \frac{1}{2} T \sum\limits_{k=0}^{n-1} [F(x_k)]^2 \,
          [\overline{u}_k^1]^2 \\
    &= \frac{1}{2} \frac{1}{T} \, h^4 \,
          \sum\limits_{k=0}^{n-1} \left[ 
          \frac{F(x_k)}{x_k} \left(
          \cos x_k - \frac{\sin x_k}{x_k} \right) \right]^2
          \left[
          \sum\limits_{i=1}^n \partial_t u(t_{i-1}) \sin ((2i-1)x_k)
          \right]^2 \\
    &\leq \frac{1}{3\pi} \frac{1}{T} \, h^4 \, 
          \sum\limits_{k=0}^{n-1} 
          \left[
          \sum\limits_{i=1}^n \partial_t u(t_{i-1}) \, x_k \, \sin ((2i-1)x_k)
          \right]^2 .
  \end{align*}
  Here, we used the following estimate, see also Figure~\ref{Fig:Fx2},
  \begin{equation} \label{AbschFx2}
    \forall x \in [0,\frac{\pi}{2}]: \quad \left[
      \frac{F(x)}{x} \left(
        \cos x - \frac{\sin x}{x} \right) \right]^2
    \leq \frac{2}{3\pi} \, x^2.
  \end{equation}
  
  \begin{figure}[ht]
    \centering
    \includegraphics[width = 0.5\textwidth]{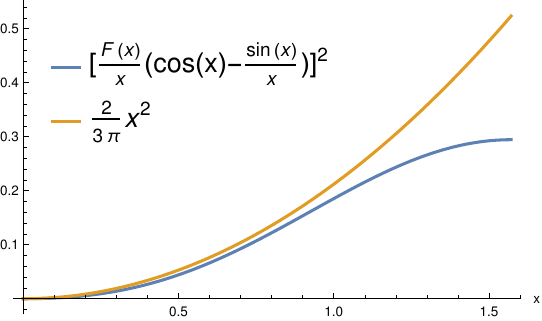}
    \caption{Justifying estimate~\eqref{AbschFx2}.} \label{Fig:Fx2}
  \end{figure}
  
  \noindent
  When using
  \[
    \frac{x}{\sin x} \leq \frac{\pi}{2} \quad \mbox{for all} \;
    x \in [0,\frac{\pi}{2}],
  \]
  we further have
  \begin{align*}
    \| w_h^1 \|^2_{L^2(0,T)}
    &\leq \frac{\pi}{12} \frac{1}{T} \, h^4 \, 
             \sum\limits_{k=0}^{n-1} \left[
             \sum\limits_{i=1}^n \partial_t u(t_{i-1}) \, \sin x_k \, \sin ((2i-1)x_k)
             \right]^2 \\
    &= \frac{\pi}{48} \frac{1}{T} \, h^4 \, 
             \sum\limits_{k=0}^{n-1} \left[
          \sum\limits_{i=1}^n \partial_t u(t_{i-1}) \Big[
          \cos (2(i-1)x_k) - \cos (2ix_k) \Big]
          \right]^2 .
  \end{align*}
  We write
  \[
    \sum\limits_{i=1}^n \partial_t u(t_{i-1}) \cos (2(i-1)x_k) =
    \partial_t u(0) + 
    \sum\limits_{i=1}^{n-1} \partial_t u(t_i) \cos (2ix_k),  
  \]
  while for $i=n$, we have
  \[
    \cos (2n x_k) = \cos \left( \frac{\pi}{2} + k\pi \right) = 0 .
  \]
  Hence, we conclude
  \[
    \sum\limits_{i=1}^n \partial_t u(t_{i-1}) \Big[
    \cos (2(i-1)x_k) - \cos (2ix_k) \Big] \, = \, \partial_t u(0) +
    \sum\limits_{i=1}^{n-1}
    \Big[ \partial_t u(t_i) - \partial_t u(t_{i-1}) \Big] \cos (2ix_k) \, ,
  \]
  and using the orthogonality
  \[
    \sum\limits_{k=0}^{n-1} \cos (2ix_k) \cos (2jx_k) = \frac{n}{2} \, \delta_{ij},
  \]
  we finally obtain
  \begin{align*}
     \sum\limits_{k=0}^{n-1}& \left[
        \sum\limits_{i=1}^n \partial_t u(t_{i-1}) \Big[
        \cos (2(i-1)x_k) - \cos (2ix_k) \Big] \right]^2 \\
     =& \sum\limits_{k=0}^{n-1} \left[ \partial_t u(0) + \sum\limits_{i=1}^{n-1}
        \Big[ \partial_t u(t_i)-\partial_t u(t_{i-1})\Big] \cos (2ix_k) \right]^2 \\
     =& \sum\limits_{k=0}^{n-1} [\partial_t u(0)]^2 +
        2 \, \partial_t u(0) \sum\limits_{i=1}^{n-1} \Big[ \partial_t u(t_i) - \partial_t u(t_{i-1}) \Big]
        \sum\limits_{k=0}^{n-1} \cos ( 2i x_k) \\
     &+ \sum\limits_{i=1}^{n-1} \sum\limits_{j=1}^{n-1}
        [\partial_t u(t_i)-\partial_t u(t_{i-1})][\partial_t u(t_j)-\partial_t u(t_{j-1})]
        \sum\limits_{k=0}^{n-1} \cos (2ix_k) \cos (2jx_k) \\
     =&  n \, [\partial_t u(0)]^2 + \frac{n}{2} \,
        \sum\limits_{i=1}^{n-1} [\partial_t u(t_i)-\partial_t u(t_{i-1})]^2 \, \leq \,
        n \, [\partial_t u(0)]^2 +
        \frac{n}{2} \, \sum\limits_{i=1}^n \left[
        \int_{t_{i-1}}^{t_i} \partial_{tt} u(s) \, ds \right]^2 \\
     \leq& \, n \, [\partial_t u(0)]^2 +
        \frac{n}{2} \sum\limits_{i=1}^n
        \int_{t_{i-1}}^{t_i} [\partial_{tt} u(s)]^2 \, ds
        \int_{t_{i-1}}^{t_i} 1^2 \, ds \, = \, n \, [\partial_t u(0)]^2 +
        \frac{nh}{2} \sum\limits_{i=1}^n
        \int_{t_{i-1}}^{t_i} [\partial_{tt} u(s)]^2 \, ds \\
     =& \, h^{-1} \, T  \, [\partial_t u(0)]^2
        + \frac{1}{2} \, T \, \int_0^T [\partial_{tt} u(s)]^2 \, ds \, .
  \end{align*}
  This concludes the proof.}

\noindent
With the above results, we are in position to give an estimate for
$w_h = Q_h {\mathcal{H}}_T^{-1}(u-Q_hu)$. When combining \eqref{Splitting}
with \eqref{Norm u2} and \eqref{Estimate wh10}, this gives
\begin{equation}\label{Bound wh H20}
  \| Q_h {\mathcal{H}}_T^{-1}(u-Q_hu) \|_{L^2(0,T)}
  \, \leq \, c \, h^2 \, \| \partial_{tt} u \|_{L^2(0,T)} \, ,
\end{equation}
for $u \in H^2(0,T)$ with $\partial_t u(0)=0$. For $u \in H^2(0,T)$ with $\partial_t u(0) \neq 0$,
we have to use \eqref{Estimate wh1} to conclude
  \begin{equation}\label{Bound wh H2}
  \| Q_h {\mathcal{H}}_T^{-1}(u-Q_hu) \|_{L^2(0,T)}
  \, \leq \, c \, \Big[ h^4 \, \| \partial_{tt} u \|^2_{L^2(0,T)} + h^3 \,
  [\partial_t u(0)]^2 \Big]^{1/2} \, .
\end{equation}
When using a space interpolation argument, we formulate these
results in a more general way.

\begin{lemma}
  Let $u \in H^s_{0,}(0,T) := [H^2_{0,}(0,T);H^1_{0,}(0,T)]_s$ with interpolation norm $\| \cdot \|_{H^s_{0,}(0,T)}$ for some
  $s \in [1,2]$, where the Sobolev space $H^2_{0,}(0,T) := \{ v \in H^2(0,T) : v(0)=\partial_t v(0)=0 \}$ is endowed with the Hilbertian norm $\| \cdot \|_{H^2_{0,}(0,T)} := \| \partial_{tt}(\cdot)\|_{L^2(0,T)}$.
  Then, the estimate
  \begin{equation}\label{Bound wh}
    \| Q_h {\mathcal{H}}_T^{-1}(u-Q_hu) \|_{L^2(0,T)}
    \, \leq \, c \, h^s \, \|u\|_{H^s_{0,}(0,T)}
  \end{equation}
  holds true.
\end{lemma}

\proof{For $s=1$ the assertion is \eqref{Estimate cB H1}, while for
  $s=2$, the assertion is \eqref{Bound wh H20}. Then, for $s \in (1,2)$,
  the assertion follows from a space interpolation argument.}

\noindent
In order to verify the estimates \eqref{Bound wh H20} and
\eqref{Bound wh H2}, we first consider $u(t) = t^3-10t^2$,
$t\in (0,T)$, $T=2$, with
$u \in H^2_{0,}(0,T)$, i.e., \eqref{Bound wh H20} implies second-order
convergence, see Table \ref{Table Bound wh} for the numerical results.
As a second example, we consider $u(t)=t^3-10t$ with $\partial_t u(0)=-10$,
where \eqref{Bound wh H2} implies a reduced order of convergence
$h^{3/2}$, as observed in the numerical example, see also
Table~\ref{Table Bound wh}.

\begin{table}[h]
\centering
    \begin{tabular}{rcccccc}
      \hline
      && \multicolumn{2}{c}{$u(t)=t^3-10t^2$}
      && \multicolumn{2}{c}{$u(t)=t^3-10t$} \\   
      $n$ && $\| Q_h {\mathcal{H}}_T^{-1}(u-Q_hu) \|_{L^2(0,T)}$
      & eoc && $\| Q_h {\mathcal{H}}_T^{-1}(u-Q_hu)\|_{L^2(0,T)}$ & eoc \\
      \hline
      2   && 1.86270658 &      && 1.80230151 &      \\
      4   && 0.45239154 & 2.00 && 0.63072107 & 1.50 \\
      8   && 0.11030583 & 2.00 && 0.21675270 & 1.50 \\
      16  && 0.02713342 & 2.00 && 0.07528792 & 1.50 \\
      32  && 0.00671882 & 2.00 && 0.02637954 & 1.50 \\
      64  && 0.00167065 & 2.00 && 0.00928626 & 1.50 \\
      128 && 0.00041642 & 2.00 && 0.00327641 & 1.50 \\
      \hline
  \end{tabular}
\caption{Numerical results for 
  $\| Q_h {\mathcal{H}}_T^{-1}(u-Q_hu)\|_{L^2(0,T)}$ with $T=2$ in the case
  $u(t)=t^3-10t^2$ with $\partial_t u(0)=0$, and for
  $u(t)=t^3-10t$ with $\partial_t u(0)\neq 0$.}
  \label{Table Bound wh}
\end{table}

\noindent
When summarizing all previous results, we state the main result
of this section.

\begin{theorem}
  Assume $u \in H^s(0,T)$ for some $s \in [0,1]$.
  Let $u_h \in S_h^0(0,T)$ be the unique solution
  of the variational formulation \eqref{Hilbert L2 Projektion}. Then, the error estimate
  \begin{equation}\label{final estimate}
    \| u - u_h \|_{L^2(0,T)} \leq
    c \, \left[ h^s \, \|u\|_{H^s(0,T)} +
      \frac{\| Q_h {\mathcal{H}}_T^{-1}(u-Q_hu)\|_{L^2(0,T)}}
      {c_S(u_h-Q_hu)} \right]
  \end{equation}
  holds true $s \in [0,1]$.
  For $u \in H^2(0,T)$ with $\partial_t u(0)=0$, we have optimal convergence, i.e.,
  \[
    \| u - u_h \|_{L^2(0,T)} \leq
    c \, h \, \| u \|_{H^2(0,T)}.
  \]
\end{theorem}

\noindent
With the results of this section, we are in the position to review the numerical examples
of Section \ref{Section:Projection}.
First, we consider the function $u(t) = \sin (\frac{\pi}{4}t)$ for
$t \in (0,2)$. We obviously have $\partial_t u(0) \neq 0$, and hence
\eqref{Bound wh H2} applies, i.e.,
$\| Q_h {\mathcal{H}}_T^{-1}(u-Q_hu)\|_{L^2(0,T)}={\mathcal{O}}(h^{3/2})$.
On the other hand, we also compute
$c_S(u_h-Q_hu) = {\mathcal{O}}(h^{1/2})$, see Table \ref{cScBreg}.
In this case, \eqref{final estimate} gives
$\| u - u_h \|_{L^2(0,T)} = {\mathcal{O}}(h)$
as already observed in Table~\ref{Table Fehler regular}.

\begin{table}[ht]
\centering
\begin{tabular}{rcccc}
\hline
  $n$ & $c_S(u_h-Q_hu)$ & eoc & $\|Q_h\mathcal{H}^{-1}_T(u_h-Q_hu))\|_{L^2(0,T)}$
  & eoc \\ \hline
   2 & 4.290 --1 &      & 1.411 --1 &      \\
   4 & 3.437 --1 & 0.32 & 4.922 --2 & 1.52 \\
   8 & 2.432 --1 & 0.50 & 1.691 --2 & 1.54 \\
  16 & 1.700 --1 & 0.52 & 5.888 --3 & 1.52 \\
  32 & 1.193 --1 & 0.51 & 2.067 --3 & 1.51 \\
  64 & 8.403 --2 & 0.51 & 7.284 --4 & 1.50 \\
 128 & 5.931 --2 & 0.50 & 2.572 --4 & 1.50 \\
 256 & 4.190 --2 & 0.50 & 9.086 --5 & 1.50 \\
 512 & 2.961 --2 & 0.50 & 3.211 --5 & 1.50 \\
1024 & 2.093 --2 & 0.50 & 1.135 --5 & 1.50 \\
2048 & 1.480 --2 & 0.50 & 4.013 --6 & 1.50 \\
\hline
\end{tabular}
\caption{Values of $c_S(u_h-Q_hu)$ and
  $\|Q_h\mathcal{H}^{-1}_T(u_h-Q_hu))\|_{L^2(0,T)}$
in the case of a regular function $u(t)=\sin (\frac{\pi}{4}t)$, $t\in (0,2)$.}
\label{cScBreg}
\end{table}

The second example was the singular function $u(t)=t^{2/3}$, $t \in (0,2)$,
i.e., $u \in H^s_{0,}(0,T)$ for $s < \frac{7}{6}$.
In this case, \eqref{Bound wh} gives
$\| Q_h {\mathcal{H}}_T^{-1}(u-Q_hu)\|_{L^2(0,T)}={\mathcal{O}}(h^{7/6})$.
The numerical results as shown in Table \ref{cSvBsing0} indicate
$c_S(u_h-Q_hu) = {\mathcal{O}}(h^{1/2})$, and therefore,
\eqref{final estimate} implies
$\| u - u_h \|_{L^2(0,T)} = {\mathcal{O}}(h^{2/3})$, as already observed
in Table \ref{Table Fehler singular 0}.

\begin{table}[ht]
\centering
\begin{tabular}{rcccc}
\hline
$n$ & $c_S(u_h-Q_hu)$ & eoc & $\|Q_h\mathcal{H}^{-1}_T(u_h-Q_hu))\|_{L^2(0,T)}$
                              & eoc \\ \hline
   2 & 4.467 --1 &      & 1.629 --1 &      \\
   4 & 2.990 --1 & 0.58 & 7.255 --2 & 1.17 \\
   8 & 2.055 --1 & 0.54 & 3.233 --2 & 1.17 \\
  16 & 1.433 --1 & 0.52 & 1.440 --2 & 1.17 \\
  32 & 1.006 --1 & 0.51 & 6.415 --3 & 1.17 \\
  64 & 7.090 --2 & 0.51 & 2.858 --3 & 1.17 \\
 128 & 5.005 --2 & 0.50 & 1.273 --3 & 1.17 \\
 256 & 3.536 --2 & 0.50 & 5.670 --4 & 1.17 \\
 512 & 2.499 --2 & 0.50 & 2.526 --4 & 1.17 \\
1024 & 1.767 --2 & 0.50 & 1.125 --4 & 1.17 \\
2048 & 1.249 --2 & 0.50 & 5.012 --5 & 1.17 \\
\hline
\end{tabular}
\caption{Values of $c_S(u_h-Q_hu)$ and
  $\|Q_h\mathcal{H}^{-1}_T(u_h-Q_hu))\|_{L^2(0,T)}$
  in the case of the function $u(t)=t^{2/3}$, $t\in (0,2)$, with a singularity
    at $t=0$.}\label{cSvBsing0}
\end{table}

Finally, we consider the function $u(t)=t(2-t)^{2/3}$, $t \in (0,2)$,
with a singularity at the terminate time $T=2$. Again we have
$\| Q_h {\mathcal{H}}_T^{-1}(u-Q_hu)\|_{L^2(0,T)}={\mathcal{O}}(h^{7/6})$,
but in this case we observe, at least asymptotically,
$c_S(u_h-Q_hu) = {\mathcal{O}}(h^{1/6})$, see Table \ref{cScBsingT}.
With these results,
\eqref{final estimate} implies
$\| u - u_h \|_{L^2(0,T)} = {\mathcal{O}}(h)$, as already observed
in Table \ref{Table Fehler singular T}.

\begin{table}[ht]
\centering
\begin{tabular}{rcccc}
\hline
$n$ & $c_S(u_h-Q_hu)$ & eoc & $\|Q_h\mathcal{H}^{-1}_T(u_h-Q_hu))\|_{L^2(0,T)}$& eoc \\ \hline
   2 & 4.128 --1 &      & 3.040 --1 &      \\
   4 & 3.391 --1 & 0.28 & 1.104 --1 & 1.46 \\
   8 & 2.599 --1 & 0.38 & 3.968 --2 & 1.48 \\
  16 & 1.952 --1 & 0.41 & 1.442 --2 & 1.46 \\
  32 & 1.476 --1 & 0.40 & 5.353 --3 & 1.43 \\
  64 & 1.138 --1 & 0.37 & 2.042 --3 & 1.39 \\
 128 & 9.009 --2 & 0.34 & 8.022 --4 & 1.35 \\
 256 & 7.326 --2 & 0.30 & 3.247 --4 & 1.31 \\
 512 & 6.106 --2 & 0.26 & 1.349 --4 & 1.27 \\
1024 & 5.193 --2 & 0.23 & 5.724 --5 & 1.24 \\
2048 & 4.485 --2 & 0.21 & 2.468 --5 & 1.21 \\
\hline
\end{tabular}
\caption{Values of $c_S(u_h-Q_hu)$ and
  $\|Q_h\mathcal{H}^{-1}_T(u_h-Q_hu))\|_{L^2(0,T)}$
  in the case of the function $u(t)=t (2-t)^{2/3}$, $t\in (0,2)$,
  with a singularity
    at $T=2$.}\label{cScBsingT}
\end{table}

\section{Conclusions} \label{Section:Con}
In this paper, we have given a complete numerical analysis to prove a
discrete inf-sup stability condition for the modified Hilbert transformation
${\mathcal{H}}_T$. While the stability constant is mesh dependent, related
error estimates are still optimal, in most cases. We restrict our
theoretical considerations to the case of piecewise constant basis
functions, however, this approach can be extended to higher-order basis functions
as well, as it is confirmed by numerical results for piecewise linear
and second-order basis functions.

These results are of utmost importance in the numerical analysis of
space-time finite element methods to analyze 
discrete inf-sup stability conditions and related error estimates
for evolution equations, for both parabolic and hyperbolic problems.
In particular, these results will provide
the stability and error analysis for a space-time finite element method
for the wave equation which is unconditionally stable. While numerical
results were already given in \cite{LoescherSteinbachZank:2022},
its numerical analysis will be given in a forthcoming paper
\cite{LoescherSteinbachZank:Wave}. Further, in the parabolic case, the 
numerical analysis in \cite{HarbrechtSchwabZank:WaermeWavlets} of a sparse grid 
approach based on wavelets also benefits from the results presented here.

\section*{Declarations}
{\bf Conflict of interest:} The authors declared that they have
no conflict of interest. \\[1mm]
{\bf Data availability:}
Data will be made available on request.

\section*{Acknowledgments}
Part of the work was done when the third author was a NAWI Graz
PostDoc Fellow at the Institute of Applied Mathematics, TU Graz.
The authors acknowledge NAWI Graz for the financial support.

\section*{Appendix}
In this appendix, we provide some of the technical computations which were
frequently used within this paper. Let $n \in \N$ be the number of finite 
elements. First, we recall the definition
\eqref{Def xk} of $x_k$. For $k=0,\ldots,2n-1$ and $\mu \in {\mathbb{N}}_0$,
we consider the shift
\[
  x_{k+2\mu n} =
  \left( \frac{\pi}{2} + (k+2\mu n)\pi \right) \frac{1}{2n} =
  x_k + \mu \pi .
\]
Then, the evaluation of certain trigonometric functions gives
\begin{align*}
  \cos (x_{k+2\mu n})
  &=\cos (x_k+\mu \pi) \\
  &=\cos (x_k) \,  \cos (\mu \pi) - \sin x_k \, \sin (\mu \pi)
        \, = \, (-1)^\mu \cos (x_k), \\
  \sin (x_{k + 2\mu n})
  &=\sin (x_k + \mu \pi) \\
  &=\sin x_k \, \cos (\mu \pi) + \cos (x_k) \, \sin (\mu \pi) \, = \,
        (-1)^\mu \sin x_k, \\
  \sin ((2j-1) x_{k+2\mu n})
  &=\sin [(2j-1)(x_k+\mu \pi)] \\
  &=\sin [(2j-1)x_k] \cos [(2j-1)\mu \pi] +
        \cos [(2j-1)x_k] \sin [(2j-1)\mu \pi] \\
  &=(-1)^\mu \sin [(2j-1)x_k] , \\
  \cos ((2j-1) x_{k+2\mu n})
  &=\cos [(2j-1)(x_k+\mu \pi)] \\
  &=\cos [(2j-1)x_k] \cos [(2j-1)\mu \pi] -
        \sin [(2j-1)x_k] \sin [(2j-1)\mu \pi] \\
  &=(-1)^\mu \cos [(2j-1)x_k] .
\end{align*}
Next, we consider, for $k=0,\ldots,n-1$,
\[
  x_{2n-1-k} = \Big( \frac{\pi}{2}+(2n-1-k) \pi \Big) \frac{1}{2n} =
  \pi - x_k,
\]
and we compute
\begin{align*}
  \cos (x_{2n-1-k})
  &=\cos (\pi-x_k) \, = \, - \cos x_k, \\
  \sin (x_{2n-1-k})
  &=\sin (\pi-x_k) \, = \, \sin x_k, \\
  \sin ((2j-1)x_{2n-1-k})
  &=\sin [(2j-1)(\pi-x_k)] \\
  &=\sin [(2j-1) \pi] \, \cos [(2j-1) x_k] -
        \cos [(2j-1) \pi] \, \sin [(2j-1) x_k] \\
  &=\sin [(2j-1) x_k] , \\
  \cos ((2j-1)x_{2n-1-k})
  &=\cos [(2j-1)(\pi-x_k)] \\
  &=\cos [(2j-1) \pi] \, \cos [(2j-1) x_k] +
        \sin [(2j-1) \pi] \, \sin [(2j-1) x_k] \\
  &=- \cos [(2j-1) x_k] .
\end{align*}
In order to prove \eqref{sum cos2}, let us recall
\cite[Equation~1.342,~4.]{GradshteynRyzhik2015},
\[
  \sum\limits_{k=1}^n \cos ((2k-1) x) =
  \frac{1}{2} \, \sin (2nx) \, \frac{1}{\cos x},
\]
which yields
\[
  \sum_{i=1}^n \cos \left( \Big(\frac{\pi}{2}+k\pi\Big)
    \frac{2i-1}{n} \right) = \frac{1}{2} \,
  \sin \Big((2k+1)\pi\Big) 
  \frac{1}{\cos \left( \Big(\frac{\pi}{2}+k\pi\Big)\frac{1}{n} \right)}
  \, = \, 0 .
\]
For the computation of the Fourier coefficients \eqref{Def uk},
we consider
\begin{align*}
  \overline{u}_k
  &=\frac{2}{T} \int_0^T u_h(t) \,
        \cos \left( \Big( \frac{\pi}{2} + k \pi \Big) \frac{t}{T} \right)
        \, dt \\
  &=  \frac{2}{T} \sum\limits_{i=1}^n u_i \int_{t_{i-1}}^{t_i} 
        \cos \left( \Big( \frac{\pi}{2} + k \pi \Big) \frac{t}{T} \right)
        \, dt \\
  &=\frac{2}{T} \sum\limits_{i=1}^n u_i
        \frac{T}{\frac{\pi}{2}+k\pi} \left[
        \sin \left( \Big( \frac{\pi}{2} + k \pi \Big) \frac{t}{T} \right)
        \right]_{t_{i-1}}^{t_i} \\
  &=\frac{2}{\frac{\pi}{2}+k\pi} \sum\limits_{i=1}^n u_i
        \left[
        \sin \left( \Big( \frac{\pi}{2} + k \pi \Big) \frac{t_i}{T} \right)
        -
        \sin \left( \Big( \frac{\pi}{2} + k \pi \Big) \frac{t_{i-1}}{T} \right)
        \right] \\
\end{align*}
and further,
\begin{align*}
  \overline{u}_k
  &=\frac{4}{\frac{\pi}{2}+k\pi} \sum\limits_{i=1}^n u_i
        \cos \left( \Big( \frac{\pi}{2} + k \pi \Big)
        \frac{t_i+t_{i-1}}{2T} \right)
        \sin \left( \Big( \frac{\pi}{2} + k \pi \Big)
        \frac{t_i-t_{i-1}}{2T} \right) \\
  &=4 \frac{\sin \left( \Big( \frac{\pi}{2} + k \pi \Big) \frac{h}{2T}
        \right)}{\frac{\pi}{2}+k\pi }
        \sum\limits_{i=1}^n u_i
        \cos \left( \Big( \frac{\pi}{2} + k \pi \Big)
        \frac{2i-1}{2n} \right) \\
  &=\frac{2}{n} \,
        \frac{\sin \left( \Big( \frac{\pi}{2} + k \pi \Big) \frac{1}{2n}
        \right)}{(\frac{\pi}{2}+k\pi)\frac{1}{2n}}
        \sum\limits_{i=1}^n u_i
        \cos \left( \Big( \frac{\pi}{2} + k \pi \Big)
        \frac{2i-1}{2n} \right) \\
  &=\frac{2}{n} \, \frac{\sin x_k}{x_k} \, 
        \sum\limits_{i=1}^n u_i \cos ((2i-1)x_k)  .
\end{align*}
In particular for $k=0,\ldots,2n-1$ and $\mu \in {\mathbb{N}}_0$, we then
conclude \eqref{uk+2mn} from
\begin{align*}
  \overline{u}_{k+2\mu n}
  &=\frac{2}{n} \, \frac{\sin(x_{k+2\mu n})}{x_{k+2\mu n}}
        \, \sum\limits_{i=1}^n u_i \, \cos ((2i-1)x_{k+2\mu n}) \\
  &=\frac{2}{n} \, \frac{(-1)^\mu \sin x_k}{x_{k+2\mu n}}
        \, \sum\limits_{i=1}^n u_i \, (-1)^\mu \cos ((2i-1)x_k) \\
  &=\frac{x_k}{x_{k+2\mu n}} \, \frac{2}{n} \, \frac{\sin x_k}{x_k}
        \, \sum\limits_{i=1}^n u_i \, \cos ((2i-1)x_k) \, = \,
        \frac{2k+1}{2k+1+4\mu n} \, \overline{u}_k .
\end{align*}
Next, for $k=0,1,\ldots,n-1$, we calculate that
\begin{align*}
  \overline{u}_{2n-1-k}
  &=\frac{2}{n} \, \frac{\sin(x_{2n-1-k})}{x_{2n-1-k}}
        \, \sum\limits_{i=1}^n u_i \, \cos ((2i-1) x_{2n-1-k}) \\
  &=- \frac{2}{n} \, \frac{\sin x_k}{x_{2n-1-k}}
        \, \sum\limits_{i=1}^n u_i \, \cos ((2i-1) x_k) \\
  &=- \frac{x_k}{x_{2n-1-k}} \,
        \frac{2}{n} \, \frac{\sin x_k}{x_k}
        \, \sum\limits_{i=1}^n u_i \, \cos ((2i-1) x_k)  \, = \,
        - \frac{2k+1}{4n-1-2k} \, \overline{u}_k^1 ,
\end{align*}
i.e., \eqref{uk+n} follows.

For the determination of the coefficients \eqref{Berechnung wj}, we compute that
\begin{align*}
  w_j &=\frac{1}{h} \int_{t_{j-1}}^{t_j}
            ({\mathcal{H}}_T^{-1} u_h)(t) \, dt \\
      &=\frac{1}{h} \sum\limits_{k=0}^\infty \overline{u}_k
            \int_{t_{j-1}}^{t_j}
            \sin \left( \Big( \frac{\pi}{2} + k \pi \Big) \frac{t}{T} \right)
            \, dt \\
      &=\frac{1}{h} \sum\limits_{k=0}^\infty \overline{u}_k
            \frac{T}{\frac{\pi}{2}+k\pi} \left[
            -  \cos \left( \Big(
            \frac{\pi}{2} + k \pi \Big) \frac{t}{T} \right)
            \right]_{t_{j-1}}^{t_j} \\
      &=\frac{1}{h} 
            \sum\limits_{k=0}^\infty \overline{u}_k
            \frac{T}{\frac{\pi}{2}+k\pi}
            \left[
            \cos \left( \Big(
            \frac{\pi}{2} + k \pi \Big) \frac{t_{j-1}}{T} \right)
            -  \cos \left( \Big(
            \frac{\pi}{2} + k \pi \Big) \frac{t_j}{T} \right)
            \right]
\end{align*}
and further,
\begin{align*}
  w_j &=\frac{1}{h} 
            \sum\limits_{k=0}^\infty \overline{u}_k
            \frac{2T}{\frac{\pi}{2}+k\pi}
            \sin \left( \Big(
            \frac{\pi}{2} + k \pi \Big) \frac{t_j+t_{j-1}}{2T} \right)
            \sin \left( \Big(
            \frac{\pi}{2} + k \pi \Big) \frac{t_j-t_{j-1}}{2T} \right) \\
      &=\frac{1}{h} 
            \sum\limits_{k=0}^\infty \overline{u}_k
            \frac{2T}{\frac{\pi}{2}+k\pi}
            \sin \left( \Big(
            \frac{\pi}{2} + k \pi \Big) \frac{2j-1}{2n} \right)
            \sin \left( \Big(
            \frac{\pi}{2} + k \pi \Big) \frac{h}{2T} \right) \\
      &=\sum\limits_{k=0}^\infty \overline{u}_k
            \frac{\sin \left( \Big(
            \frac{\pi}{2} + k \pi \Big) \frac{h}{2T} \right)}
            {(\frac{\pi}{2}+k\pi)\frac{h}{2T}}
            \sin \left( \Big(
            \frac{\pi}{2} + k \pi \Big) \frac{2j-1}{2n} \right) \\
      &=\sum\limits_{k=0}^\infty \overline{u}_k \,
            \frac{\sin x_k}{x_k} \,
            \sin ((2j-1)x_k) .
\end{align*}

To prove \eqref{wj+2mn}, we have
\begin{align*}
  w_j &=\sum\limits_{k=0}^\infty \overline{u}_k \, \frac{\sin x_k}{x_k}
            \, \sin((2j-1)x_k) \\
    &=\sum\limits_{k=0}^{2n-1} \sum\limits_{\mu=0}^\infty
          \overline{u}_{k+2\mu n} \, \frac{\sin(x_{k+2\mu n})}{x_{k+2\mu n}} \,
          \sin ((2j-1)x_{k+2\mu n}) \\
    &=\sum\limits_{k=0}^{2n-1} \sum\limits_{\mu=0}^\infty
          \overline{u}_{k+2\mu n} \, \frac{(-1)^\mu\sin x_k}{x_{k+2\mu n}} \,
          (-1)^\mu \sin ((2j-1)x_k) \\
    &=\sum\limits_{k=0}^{2n-1} \overline{u}_k
          \left(  \sum\limits_{\mu=0}^\infty
          \frac{(2k+1)^2}{(2k+1+4\mu n)^2} \right)
          \, \frac{\sin x_k}{x_k} \,
            \sin ((2j-1)x_k) \\
    &=\sum\limits_{k=0}^{2n-1} \overline{u}_k \, \gamma(k,n)
          \, \frac{\sin x_k}{x_k} \,
          \sin ((2j-1)x_k)
  \end{align*}
  when using \eqref{uk+2mn}.
  
  Next, we use \eqref{uk+n} to conclude from
  \eqref{wj+2mn} that
  \begin{align*}
    w_j =&\sum\limits_{k=0}^{2n-1} \overline{u}_k \, \gamma(k,n) \,
              \frac{\sin x_k}{x_k} \, \sin ((2j-1) x_k) \\
        =&\sum\limits_{k=0}^{n-1} \left[
              \overline{u}_k \, \gamma(k,n) \, \frac{\sin x_k}{x_k}
              \sin ((2j-1)x_k) 
              \right. \\
        & \left.
              + \overline{u}_{2n-1-k} \, \gamma(2n-1-k,n) \,
       \frac{\sin(x_{2n-1-k})}{x_{2n-1-k}} \,
       \sin ((2j-1) x_{2n-1-k}) 
              \right] \\
       =&\sum\limits_{k=0}^{n-1} \overline{u}_k \Big[
              \gamma(k,n) 
              - \gamma(2n-1-k,n)
              \frac{(2k+1)^2}{(4n-1-2k)^2}
              \Big]
              \frac{\sin x_k}{x_k} \,
              \sin ((2j-1) x_k),
  \end{align*}
  which yields \eqref{wj+n}.

\end{document}